\documentclass[lettersize,journal]{IEEEtran}
\usepackage[utf8]{inputenc}
\usepackage[T1]{fontenc}

\usepackage{enumitem}
\usepackage{amsmath, amssymb}
\usepackage{semantic}
\usepackage{nicematrix}
\usepackage[caption=false,font=normalsize,labelfont=sf,textfont=sf]{subfig}
\usepackage{arydshln}
\usepackage{algpseudocode}
\usepackage{algorithm}
\usepackage{graphicx}
\usepackage{comment}
\usepackage{balance}
\usepackage{optidef}
\usepackage{mathrsfs}
\usepackage{hyperref}
\usepackage{cite}
\usepackage{varwidth}
\usepackage{tikz,xcolor,hyperref}
\definecolor{lime}{HTML}{A6CE39}
\DeclareRobustCommand{\orcidicon}{
	\begin{tikzpicture}
	\draw[lime, fill=lime] (0,0) 
	circle [radius=0.16] 
	node[white] {{\fontfamily{qag}\selectfont \tiny ID}};
	\draw[white, fill=white] (-0.0625,0.095) 
	circle [radius=0.007];
	\end{tikzpicture}
	\hspace{-2mm}
}
\foreach \x in {A, ..., Z}{\expandafter\xdef\csname orcid\x\endcsname{\noexpand\href{https://orcid.org/\csname orcidauthor\x\endcsname}
			{\noexpand\orcidicon}}
}
\usepackage[compact]{titlesec}         
\titlespacing{\section}{8pt}{8pt}{8pt} 
\AtBeginDocument{
  \setlength\abovedisplayskip{2pt}
  \setlength\belowdisplayskip{2pt}} 
\newtheorem{theorem}{Theorem}
\newtheorem{proposition}{Proposition}
\newtheorem{lemma}{Lemma}

\newtheorem{remark}{Remark}
\newtheorem{assumption}{Assumption}

\newtheorem{corollary}{Corollary}
\newtheorem{bound}{Bound}

\newcommand{\R}{\mathbb{R}}
\newcommand{\Rn}{\mathbb{R}^n}
\newcommand{\Rnn}{\mathbb{R}^{n \times n}}
\newcommand{\qed}{\hfill $\square$}

\newcommand{\bac}[1]{\textcolor{blue}{#1}}

\def\BibTeX{{\rm B\kern-.05em{\sc i\kern-.025em b}\kern-.08em
    T\kern-.1667em\lower.7ex\hbox{E}\kern-.125emX}}
\begin{document}
\title{Online Distributed Algorithm for Optimal Power Flow problem with Regret Analysis}
\author{Sushobhan Chatterjee\orcidA{} and Rachel Kalpana Kalaimani\orcidB{}
\thanks{*This work has been partially supported by DST-INSPIRE Faculty Grant, Department of Science and Technology (DST), Govt. of India (ELE/16-17/333/DSTX/RACH)}
\thanks{The authors are with the Department of Electrical Engineering, Indian Institute of Technology Madras, Chennai, India - 600036. \texttt{(email: chatterjeesushobhan8@gmail.com; rachel@ee.iitm.ac.in)}.}%
}

\maketitle

\begin{abstract}
We investigate the distributed DC-Optimal Power Flow (DC-OPF) problem for a dynamic and uncertain environment. The unpredictable supply of renewable resources and varying prices of the electricity market are a few factors responsible for the uncertainty. 
We propose to address this problem using the framework of online convex optimization, where the cost functions are not known apriori because of the uncertainty and are revealed only incrementally over time. We also consider a  distributed setting, where each agent (generators and loads) in the power network is only privy to their own local objectives and constraints but can communicate with their neighbours.
 A distributed online algorithm is
proposed based on the modified primal-dual approach. The performance of the online algorithm is evaluated using the regret (static) function, which is the difference between the actual cost incurred by employing the proposed algorithm and the optimal fixed decision in hindsight. Since we deal with a constrained optimization problem, analogous to the notion of regret the accumulation of the constraint violation is also calculated at each step. We establish a sub-linear bound on the static regret and constraint violation under suitable assumptions on step-size and cost function. Finally, we use the standard IEEE-14 bus system to demonstrate the performance of our algorithm. 
\end{abstract}


\section{Introduction} \label{intro}
\IEEEPARstart{I}{n Recent} times, with the rapid advancements in renewable energy resources and electric vehicles, there has been a colossal 
 transformation in the 
 electrical power grid; a dynamically changing and uncertain environment is increasingly inundating the picture \cite{dumlao2022impact, smith2022effect}. 
 This is due to the unpredictable nature of the renewable
 generating resources viz. wind \cite{athari2017impacts} and solar \cite{nwaigwe2019overview}, the impact of harmonics due to the presence of power conditioners \cite{khan2020stability}, volatility of the electricity prices due to real-time spot market bids \cite{roozbehani2012volatility}, net metering \cite{abdin2018electricity}, and so forth. 
 With this inception,
 the way of tackling classical problems like economic load dispatch \cite{bai2017distributed}, demand response \cite{ye2016game}, and optimal power flow (OPF) \cite{frank2012optimal} needs a revamp.

The classical (offline) OPF problem does not take into account the above mentioned uncertainties inherent in the power grid 
on a day-to-day basis, which significantly influences the  generation cost.  Hence, future algorithms for OPF require a real-time (online) operation, so as to account for the uncertainty.
Our work investigates the OPF problem using an \textit{online learning-based optimization} framework. 
It is also being used in numerous other applications such as portfolio optimization, recommender engines, auctions, and microgrids, to name a few \cite{zinkevich2003online, wang2022designing, chadoulos2021learning}.

The basic tenet of an online problem is the dynamic and uncertain nature of its objective function and possibly constraints, which are sequentially unveiled to the decision maker. Hence, the decision maker operates with past information to take a decision at the current  time step. The consequence of the decision is then revealed to the user which enables them to take better decisions in future and hopefully converge to the optimal decision in the end \cite{li2022survey}.
The performance of an online algorithm is judged on the basis of a metric called \textit{regret}; which measures how "good" the online algorithm is in meeting objectives compared to its offline counterpart, where the environment is completely known prior.
In the case of a constrained optimization setting, 
another metric viz. \textit{violation} is employed which tracks the extent of constraint violations. An online algorithm is termed acceptable if the growth of the above two metrics is sub-linear \cite{hazan2016introduction}.

Most of the online algorithms in the literature are for "inequality-constrained" optimization problems. However a very important class of dynamic problems involving optimized resource allocation as seen in the power systems like economic load dispatch \cite{bai2017distributed}, optimal power flow \cite{mohammadi2014distributed} etc., are primarily "equality-constrained" in nature. Unlike inequality constraints, the fundamental nature of the dual variables associated with equality constraints at optimality remains \textit{unknown} which makes the regret analysis fundamentally unique, and as such very limited literature exists in this domain. The only recent work \cite{lupien2023online} which addresses this problem uses a centralized $2^{nd}$ order optimality method, however for a large-scale system like power networks this brings in a higher complexity in the form of high computational cost $O(N^3)$, saddle-convergence problems \cite{haeser2018some} and scalability issues \cite{cartis2018second}. We propose a lesser complex algorithm using $1^{st}$ order methods which properly scales up with the order of the multi-agent systems.

The recent interest in developing distributed algorithms to solve power flow problems has been gaining momentum for the past few years and stems primarily from the advent of large-scale distributed generation strategies and intelligent demand management. Hence, it becomes necessary to find ways to be able to accommodate such a massive increase in variables within a specific control area, and also improve coordination among different control areas. In such a scenario, distributed methods offer increased reliability and easy-to-implement ability, primarily on account of the low-cost systems and advanced computation algorithms \cite{li2019distributed, aybat2017distributed}. Hence we propose a distributed online algorithm for the OPF problem.



\subsection{Related Work} 

In the area of power systems, online learning is increasingly being utilized for informed decision-making under uncertainty. This uncertainty usually manifests itself through system models or objectives. Examples include social-welfare maximization using online feedback \cite{zhou2017incentive}, 
charge scheduling of electric vehicles 
\cite{vujanic2016decomposition}, online distributed energy management using microturbines  \cite{ma2016distributed},  online cost minimization of wind turbines against a spot market transaction scenario \cite{narayanaswamy2012online}, and online economic load dispatch
\cite{badiei2015online}.

While our paper handles OPF using an online optimization framework, a closely related approach in the literature is \textit{online feedback-based optimization} \cite{zhu2020distributed, bernstein2019online, gan2016online, dall2016optimal}. 
In this approach, the cost functions (either static or time-varying) are known prior, and a feedback control structure is used to gather real-time measurements to account for the imperfections in the system models. On the other hand, we employ \textit{online learning-based optimization} framework which works with past data to "learn" the adversarial or stochastic nature of the dynamic optimization problem. In comparison to the former case, here the cost functions are not available while taking a decision. Even though the former approach provides robustness against uncertainty, in essence, is different from the latter approach which deals with learning the uncertainty.  Prior works on feedback framework include online algorithms for tracking set-points \cite{zhu2020distributed}, feedback-based online algorithms for AC-OPF problem \cite{bernstein2019online}, using gradient methods for OPF with bounded sub-optimality gap \cite{gan2016online}, relaxed AC-OPF feedback algorithms \cite{dall2016optimal}, etc.


An attempt has been made in \cite{kim2014online} for investigating the centralized online learning-based setting for OPF problems without any guarantees on regret bounds. 
To the best of our knowledge, our work  is the first  on distributed online optimization for the DC-OPF problem with regret analysis yielding theoretical guarantees on the performance of the proposed algorithm. 
\subsection{Contributions and Outline}
Our primary contributions can be summarized as follows.
\begin{enumerate}
    \item We propose a distributed primal-dual-based "online convex optimization (OCO)" algorithm which learns and solves the uncertainty-laden DC optimal power flow (DC-OPF) problem for a power transmission network, with the agents communicating over an undirected graph.
    \item We provide theoretical guarantees for the convergence of our algorithm via. online performance metrics viz. \textit{static} regret and constraint violation, and establish sublinear bounds of order $O(\sqrt{T})$ for them, over the time horizon $T$.
    \item We elucidate the results through numerical examples by considering the standard \textit{IEEE-14 bus system}.

\end{enumerate}

The outline of the paper is organised as follows. 
Section \ref{prelims} introduces notational conventions and brief mathematical notions. In Section \ref{pform}, we discuss the formulation of the distributed online constrained optimal power flow problem. Our proposed algorithm required for solving the problem has been posed in Section \ref{algo}. We introduce associated assumptions on the problem and theoretically establish the sublinearity of regret and violation bounds sustained by our algorithm in Section \ref{regret}. Section \ref{numex} talks about a numerical example to validate our results. Finally, we summarize and conclude our work in Section \ref{concl}. The proofs of some of the results are reported in Appendix.   

\section{Preliminaries} \label{prelims}
In this section, we will provide a few basic notational conventions and some graph-theoretical notions that will be used throughout the paper.  

\subsection{Notations}
Let $\R$, $\R_{++}$, $\Rn$ and $\Rnn$ be the set of real numbers, strictly-positive real numbers, $n$-dimensional real vector, and real square matrix of order $n$ respectively. 
For any time-invariant vector, $p$ : $p_i$ denotes its $i^{th}$ element, and for any time-dependent vector, $x_t$ : $x_{i,t}$ denotes its  $i^{th}$ element at time instant $t$, $\overline{x}_t$ denotes its average value over all $i$ at time instant $t$. Standard Euclidean vector norm is denoted by $\Vert. \rVert$, unless stated otherwise. For a matrix, $P$ : $\lVert P \rVert_{\infty}$ denotes the absolute maximum among the row sums of $P$. For a set $\mathcal{S}$, $|\mathcal{S}|$ denote its cardinality, and $\small[\mathcal{S}\small]$ refers to the elements of the set. 
For  $h\in\mathbb{R}$, $\mathbb{P}$ denotes Euclidean projection of $h$ onto its feasible region $H$, i.e., 
\begin{equation*}
\mathbb{P}_{H} \big[h \big] = \underset{f \in H}{argmin}\ \lVert f - h \rVert_2^2    
\end{equation*}

\subsection{Convex Analysis}
A set $D \subset \mathbb{R}^q$ is said to be convex if $\exists\ x, y \in D\ s.t.$
\begin{equation*}
(1-\theta)x + \theta y \in D\ \forall\ \theta \in \small[0, 1\small]
\end{equation*}
A function $f : \mathbb{R}^q \rightarrow \mathbb{R}$ is said to be convex in set $D$, if $D$ is convex for all $\theta \in \small[0, 1\small]$, and $x, y \in D$
\begin{equation*}
f(\theta x + (1-\theta) y) \leq \theta f(x) + (1-\theta) f(y)
\end{equation*}

\subsection{Graph Theory}
Let $\mathcal{G} := (\mathcal{V},\mathcal{E})$ denote an undirected graph, where $\mathcal{V} = [v_1,v_2,...,v_n]$ represents the vertex set and $\mathcal{E} \subseteq \mathcal{V}\times \mathcal{V}$ denotes the edge set. 
An undirected edge is an unordered pair of distinct vertices $(v_i,v_j)$.
The edges of the graph are sometimes assigned weights and the graph is denoted as $\mathcal{G} := (\mathcal{V},\mathcal{E},W)$, where
 the weight associated with the edge $(v_i,v_j)$ is the entry $W_{ij} \geq 0$, in $W$.  Also, we assume $W_{ii} > 0\ \forall\ i \in \mathcal{V}$. 
For a vertex $v_i$, the set of neighbours is defined by $\mathcal{N}_i=\big\{v_j\ |\ (v_i,v_j)\in \mathcal{E}\big\}$.
An undirected graph is said to be connected if there exists a sequence of edges between any two distinct pairs of vertices $\{v_i,v_j\} \in \mathcal{V}$.

\section{Problem Formulation} \label{pform}
Our primary objective lies in solving the DC-approximated optimal power flow problem under an online optimization setting. In this context, it becomes necessary to introduce the models of the underlying system which eventually paves the way towards the formulation of the online optimization problem.

\subsection{Power Distribution Networks}
Consider a power distribution topology which can be inherently modelled as an undirected graph $\mathcal{G} := (\mathcal{V}, \mathcal{E})$, with node set $\mathcal{V} := N^{-} \cup \{0\}$; $N^{-} := \{1,2,\hdots,N-1\}$, and $\mathcal{E} \subseteq \mathcal{V} \times \mathcal{V}$ edges. For the remainder of the paper, we will refer to the set $\mathcal{V}$ as $\small[N\small]$ and hence $|\mathcal{V}| := N$. Each node $i \in \small[N \small]$, will be referred to as "bus" or "agent", and edge $(k,l) \in \mathcal{E}\ \forall\ k,l \in \small[N\small]$, as "line" or "connector". Furthermore, it is considered that each bus $i \in \small[N\small]$ is having a per unit active power injection of $p_{i,t} = p^g_{i,t} - p^l_{i,t}$; with $p^g_{i,t}$ and $p^l_{i,t}$ denoting generation and load respectively, and $\theta_{i,t}$ be the corresponding load angle of the bus, at time instant $t$. 
Let bus $0$ be taken as the slack or swing bus, and it is assumed by convention that any power mismatch in the system is taken care of by it during load flow studies.

We consider the classical DC load flow model for our study, 
\begin{equation} \label{eqn:dc}
    p_i = \overset{N}{\underset{j=1}{\Sigma}} B_{ij} (\theta_i - \theta_j)
\end{equation}
where each line $(i,j) \in \mathcal{E}$ is considered to have a susceptance of $B_{ij}$; primarily due to its simple albeit powerful approximation capabilities, thus making it closer in performance to the original non-linear models \cite{stott1971effective, bienstock2019strong}.

\subsection{Generation and Load Model}
In a power system topology, generators play the role of dispatchable sources of energy and are distributed throughout the network, while the loads act as energy sinks to model consumer demands, over time. It is an underlying assumption that all the generators are readily dispatchable so as to meet the additional demands and losses of the system. The cost curves of the generation associated with thermal and renewable systems can be modelled as a purely quadratic convex function of generated power \cite{shiltz2015integrated} i.e.,
\begin{equation} \label{eqn:objf}
    f_i(p^g_i) = a_i(p^g_i)^2 + b_i(p^g_i) + c_i,
\end{equation}
where $a_i \in \mathbb{R}_{++}$, $b_i \in \mathbb{R}$ and $c_i \in \mathbb{R}$ are the incremental, base and constant cost parameters of the generator connected at the $i^{th}$ bus, respectively.
Furthermore, utility curves of customer demands can be modelled as constant-power loads $p^l_i$, denoting active load power at $i^{th}$ bus \cite{kersting2022distribution}.
\begin{remark}
Note that an $i^{th}$ bus can generally have multiple generators and loads connected to it, and can thus be aggregated to provide a coupled value viz. $p^g_i$ and, $p^l_i$ respectively.  
\end{remark}

\subsection{Optimal Power Flow (DC-OPF) Formulation : An Online Optimization approach}
The goal of the classical optimal power flow problem is related to optimizing the decision variables associated with every generator in the power system under some practical constraints, i.e., to minimize the cost of collective generation. 
\begin{mini*}|1|
 {\substack{p^g_i}}{\overset{N}{\underset{i=1}{\Sigma}} f_i(p^g_i)}{}{} \tag{$A1$}
 \label{eqn:A1}
 \addConstraint{p^g_i = p^l_i + \theta_i\underset{j \in \mathcal{N}_i}{\Sigma}B_{ij} - \underset{j \in \mathcal{N}_i}{\Sigma}B_{ij} \theta_j} \tag{$A2$} \label{eqn:A2}
 \addConstraint{0 \leq p^g_i \leq P^g_{i,max}} \tag{$A3$} \label{eqn:A3}
 \addConstraint{\theta_0 = 0,\ \forall\ i \in \small[N\small]} \tag{$A4$} \label{eqn:A4}
\end{mini*}
Here, \eqref{eqn:A1} is the objective function $f_i : P_i \rightarrow \mathbb{R}$, as defined in \eqref{eqn:objf}. \eqref{eqn:A2} describes the power-balance equation of the network as defined in \eqref{eqn:dc} and \eqref{eqn:A3} depicts the active power generation limits of the $i^{th}$ generator. Expression \eqref{eqn:A4} puts a constraint on the load angle of the slack bus, so as to set a reference.
\begin{remark}
The OPF problem \eqref{eqn:A1}-\eqref{eqn:A4} that we are dealing with is an abridged version of the original problem as it relaxes other constraints like line congestion, unit-commitment, contingency constraints etc., like in \cite{gan2016online}. 
\end{remark}
The online version of the OPF problem that is of interest to us is not just a spatial optimization over the agents as defined in \eqref{eqn:A1}-\eqref{eqn:A4}, but also a temporal one over a finite time horizon $T$, which is formulated below with appropriate modifications.
\begin{mini*}|1|
 {\substack{p^g_{i,t}}}{\overset{N}{\underset{i=1}{\Sigma}} \overset{T}{\underset{t=1}{\Sigma}} f_{i,t}(p^g_{i,t})}{}{} \tag{$B1$}
 \label{eqn:D1}
 \addConstraint{h_{i,t}(.) = 0} \tag{$B2$} \label{eqn:D2}
 \addConstraint{0 \leq p^g_{i,t} \leq P^g_{i,max}} \tag{$B3$} \label{eqn:D3}
 \addConstraint{\theta_{0,t} = 0,\ \forall\ i \in \small[N\small],\ t \in \small[T\small]} \tag{$B4$} \label{eqn:D4}
\end{mini*}
where, $f_{i,t}(p^g_{i,t}) = a^t_i(p^g_{i,t})^2 + b^t_i(p^g_{i,t}) + c^t_i$ and $h_{i,t}(.) = p^g_{i,t} - p^l_{i,t} - \theta_{i,t}\underset{j \in \mathcal{N}_i}{\Sigma}B_{ij} + \underset{j \in \mathcal{N}_i}{\Sigma}B_{ij} \theta_{j,t}$. 

The literature on online optimization theory can be utilized to model dynamically variant yet uncertain processes as a series of time-varying objective functions as depicted in \eqref{eqn:D1}. More precisely, every agent $i \in \small[N\small]$ in the network, at each time instant $t \in \small[T\small]$ takes a decision on its collective state $(p^g_{i,t},\theta_{i,t})$ without any prior information on the current cost $f_{i,t}$. After the decision-making process is over, the cost function $f_{i,t} : P_i \rightarrow \mathbb{R}$ is revealed by the nature or an adversary, and agent $i$ incurs the cost $f_{i,t}(p^g_{i,t})$. 
\begin{remark}
In the power flow studies, any imbalance in the generation is usually taken care of by the slack bus, and as such the constraint violations i.e. $h_{i,t} \neq 0$ that are ensued at each time instant $t$ while running the online OPF algorithm are also abetted by the slack bus, thereby effectively bringing down the violation to 0, for all practical purposes. Our proposed algorithm ensures that the constraints are satisfied in the \textit{long-term}.  
\end{remark}
Note that the above formulation is in a centralized form. However, the distributed structure of the OPF problem is inherent through multi-agent coordination, as can be observed in the separable nature of the objective function and only-neighbour-dependent power balance constraint of agent $i$. 

\subsection{Constrained Online Learning over Networks}

In an online optimization framework, 
due to the sequential nature of the decision-making process, the performance of the algorithm is determined by some special metrics \cite{zinkevich2003online}. In such a context, the goal of an online optimization policy is to minimize \textit{regret}. 
For any finite time horizon $T > 0$, the \textit{regret} is defined as,
\begin{equation} \label{eqn:reg}
    \mathcal{R}_s(T) = \overset{T}{\underset{t=1}{\Sigma}} \overset{N}{\underset{i=1}{\Sigma}} f_{i,t}(p^g_{i,t})\ -\ \overset{T}{\underset{t=1}{\Sigma}} \overset{N}{\underset{i=1}{\Sigma}} f_{i,t}(p^{g*}_i) 
\end{equation}
where, $p^{g*}_i \triangleq \big[\underset{p_i^g\ \forall\ i}{argmin}\  \underset{t}{\Sigma} \underset{i}{\Sigma} f_{i,t}(p^g_{i,t})\big]_i$

The above expression depicts that the decision trajectory of an online policy is being compared against the \textit{static} optimal trajectory in hindsight. Thus, it can be inferred that the quantity \textit{regret} captures how much regret an optimization policy has when working in an online setting, contrary to an offline setting with complete information.

Intuitively, it can be observed that even though the \textit{regret} will keep on inflating as time horizon $T$ expands, as long as the inflation rate is bounded to be sublinear, i.e., $\mathcal{R}_s(T)$ = $o(T)$, the algorithm performs well on average as $\underset{T \rightarrow \infty}{lim}\frac{\mathcal{R}_s(T)}{T} \rightarrow 0$. This implies that the online policy averaged over the horizon will eventually catch up with the \textit{static} optimal strategy in hindsight.

Likewise, another goal of the online policy in a constrained optimization setting is to reduce the \textit{violation} associated with constraints. 
For any arbitrary horizon $T > 0$ it is defined as,
\begin{equation} \label{eqn:ecviol}
    \mathcal{R}_s^{ec}(T) = \Bigg\lVert \overset{T}{\underset{t=1}{\Sigma}} \overset{N}{\underset{i=1}{\Sigma}} \big[h_{i,t}(.)\big] \Bigg\rVert
\end{equation}
Similar to \textit{regret}, as long as the accumulation rate of \textit{violation} is bounded to be sublinear, i.e., $\mathcal{R}_s^{ec}(T)$ = $o(T)$, the time-averaged online policy over the horizon will ensure that the equality constraints will eventually be satisfied as $\underset{T \rightarrow \infty}{lim}\frac{\mathcal{R}^{ec}_s(T)}{T} \rightarrow 0$.

\section{Online Distributed Primal-Dual based Algorithm} \label{algo}
This section will introduce a distributed online algorithm for the OPF problem which combines the distributed approach for solving the offline DC-OPF problem given by \cite{mohammadi2014distributed} with an online multi-agent coordination framework.

In this distributed online optimization framework, at any time instant $t$, each agent $i \in \small[N\small]$ possess local knowledge associated with $f_{i,t}$ and $h_{i,t}$, and has information about its neighbours $\mathcal{N}_i$.
The augmented Lagrangian function $\mathcal{L}_t : \mathbb{R}^n \times \mathbb{R}^n \times \mathbb{R}^n \times \mathbb{R}^n \times \mathbb{R} \rightarrow \mathbb{R}$ of the problem \eqref{eqn:D1}-\eqref{eqn:D4} at any time instant $t$ is given by:
\begin{equation}\label{eqn:lang}
    \begin{aligned}
        & \mathcal{L}_t = \overset{N}{\underset{i=1}{\Sigma}} f_{i,t}(p^g_{i,t}) + \overset{N}{\underset{i=1}{\Sigma}} \mu_{i,t}^{+} (p^g_{i,t} - P^g_{i,max}) + \overset{N}{\underset{i=1}{\Sigma}} \mu_{i,t}^{-} (-p^g_{i,t})
        \\&\
        + \lambda_{0,t}\theta_{0,t} + \overset{N}{\underset{i=1}{\Sigma}} \lambda_{i,t} (p^g_{i,t} - p^l_{i,t} - \theta_{i,t}\underset{j \in \mathcal{N}_i}{\Sigma}B_{ij} + \underset{j \in \mathcal{N}_i}{\Sigma}B_{ij} \theta_{j,t})
    \end{aligned}
\end{equation}
Based on the KKT stationarity and feasibility conditions associated with the problem \eqref{eqn:D1}-\eqref{eqn:D4}, at the saddle points ($*$) of $\mathcal{L}_t$ we have,
\begin{equation}\label{eqn:e1}
\nabla_{p^g_{i,t}} \mathcal{L}_t|_{*} = \nabla f_{i,t}(p^{g*}_{i,t}) + \mu_{i,t}^{+} - \mu_{i,t}^{-} + \tilde{\lambda}^{*}_{i,t} = 0
\end{equation}
\begin{equation}\label{eqn:e2}
\nabla_{\theta_{i,t}}\mathcal{L}_t|_{*} = -\tilde{\lambda}^{*}_{i,t}\underset{j \in \mathcal{N}_i}{\Sigma}B_{ij} + \underset{j \in \mathcal{N}_i}{\Sigma}B_{ij}\tilde{\lambda}^{*}_{j,t} = 0
\end{equation}
\begin{equation}\label{eqn:e3}
\nabla_{\lambda_{i,t}}\mathcal{L}_t|_{*} = p^{g*}_{i,t} - p^l_{i,t} - \theta^{*}_{i,t}\underset{j \in \mathcal{N}_i}{\Sigma}B_{ij} + \underset{j \in \mathcal{N}_i}{\Sigma}B_{ij} \theta^{*}_{j,t} = 0 
\end{equation}
\begin{equation}\label{eqn:e4}
\nabla_{\lambda_{0,t}}\mathcal{L}_t|_{*} = \theta^{*}_{0,t} = 0
\end{equation}
\begin{equation}\label{eqn:e5}
\nabla_{\mu_{i,t}^{+}}\mathcal{L}_t|_{*} = p_{i,t}^{g*} - P_{i,max}^g \leq 0 
\end{equation}
\begin{equation}\label{eqn:e6}
\nabla_{\mu_{i,t}^{-}}\mathcal{L}_t|_{*} = -p_{i,t}^{g*} \leq 0 
\end{equation}
It is observed that the dual variable linked with the constraint \eqref{eqn:D4} vanishes due to the choice of the slack bus having no impact on the final result, hence it can be omitted from the $\mathcal{L}_t$ expression. Similarly, the dual variables $\mu^{+}_{i,t}$ and $\mu^{-}_{i,t}$ can be incorporated in the primal update steps of the algorithm, rather than being dualized. 
Using the above equations, an algorithm is provided to solve the optimization problem in \eqref{eqn:D1}- \eqref{eqn:D4} in Algorithm \ref{algo:doopt}. The intuition behind the update rules \eqref{eqn:e9}-\eqref{eqn:e10} are explained in the following remark.

\begin{algorithm} 
  \caption {Local estimation of the optimal power flow by agent $i \in \small[N\small]$ using Online Primal-Dual approach}
  \label{algo:doopt}
    \begin{algorithmic}[1]
    \State \textbf{Given :} $T \geq 1$, $W$ 
    \State \textbf{Initialize :} $t=1$, $p_{i,0}^g = 0$, $\theta_{i,0} = 0$ and $\lambda_{i,0} = 0$ $\forall\ i \in \small[N\small]$
    \State \textbf{Iterate :} $\forall\ t \in \small[T \small]$ and $i \in \small[N \small]$
        \For{$t \leq T$}
            \For{$i \leq N$}
                \State Send $\lambda_{i,t}$ to out-neighbours, and receive 
                \State $\lambda_{j,t}\ \forall\ j \in \mathcal{N}_i$ from the in-neighbours.
                \begin{equation}\label{eqn:e9} \tilde{\lambda}_{i,t} = \overset{N}{\underset{j=1}{\Sigma}} W_{ij} \lambda_{j,t}
                \end{equation}
                \State \textbf{Primal Update :}
                \begin{equation}\label{eqn:e7} p_{i,t+1}^g = \mathbb{P}_{P_i} \big[p_{i,t}^g - \delta_t \nabla_{p^g_{i,t}}\mathcal{L}_t|_{\tilde{\lambda}_{i,t}}\big]
                \end{equation}
                \begin{equation}\label{eqn:e8} \theta_{i,t+1} = \theta_{i,t} - \gamma_t \nabla_{\lambda_{i,t}}\mathcal{L}_t|_{\tilde{\lambda}_{i,t}}
                \end{equation}
                \State \textbf{Dual Update :}
                \begin{equation}\label{eqn:e10} \lambda_{i,t+1} = \tilde{\lambda}_{i,t} - \beta_t \nabla_{\theta_{i,t}}\mathcal{L}_t|_{\tilde{\lambda}_{i,t}} + \alpha_t \nabla_{\lambda_{i,t}}\mathcal{L}_t|_{\tilde{\lambda}_{i,t}}
                \end{equation}
            \EndFor
        \EndFor  
    \end{algorithmic}
\end{algorithm}
\begin{remark}
Each agent $i \in \small[N\small]$ maintains four variables with it viz. local primals $p^g_{i,t} \in P_i$, $\theta_{i,t} \in \Phi_i$, local temporary $\tilde{\lambda}_{i,t} \in \mathbb{R}$, and local dual $\lambda_{i,t} \in \mathbb{R}$. The primal update for $p^g_{i,t}$ has been adopted from the standard Arrow–Hurwicz–Uzawa algorithm \cite{arrow1958studies} as seen in \eqref{eqn:e7}.
For an \textit{uncongested} transmission line i.e., one having no power delivery limits, the Lagrange multipliers $\lambda_{i,t}$ associated with the power-balance equation will eventually come into a consensus \cite{raikar2019renewable}, which is exhibited by a convex combination update given by \eqref{eqn:e9}. The dual update rule for $\lambda_{i,t}$ involves optimality and innovation terms. Intuitively, the first term as given by \eqref{eqn:e2} enforces the inherent coupling between the Lagrange multipliers to facilitate information exchange among the agents over the network, while the second term given by \eqref{eqn:e3} ensures penalizing any violations in the generation-load balance of the $i^{th}$ bus (agent), as exhibited in \eqref{eqn:e10}.
\end{remark}

\section{Regret Analysis} \label{regret}
This section will focus on the analytical treatment of the regret and violation associated with our proposed Algorithm \ref{algo:doopt}. Some key assumptions on the problem have been presented, and a few relations proposed. Finally, all the lemmas are collated to opine the main results.

\subsection{Assumptions}
There are certain standard assumptions required to ensure that the regret and violation terms are bounded for an online optimization problem.  
\begin{assumption} \label{ass:a1}
The functions $\{f_{i,t}\}_{i \in \small[N \small],\ t \geq 1}$  are convex and  differentiable.
\end{assumption}    

\begin{assumption} \label{ass:a2}
The sets $P_i$ for $i \in \small[N \small]$ are compact and time-invariant. That is, there exists a constant $C_P$ s.t.
    \begin{equation}\label{eqn:e11}
    \lVert p^g_i \rVert \leq C_P\ \forall\ p^g_i \in P_i\ \text{and i} \in \small[N \small]
    \end{equation}
    where $C_P = \underset{i}{max} \{P_{i,max}^g\}.$
    
Similarly, for the sets, $\Phi_i$ for $i \in \small[N \small]$ being compact and time-invariant, there exists a constant $C_{\theta}$ s.t. 
    \begin{equation}\label{eqn:e13}
    \lVert \theta_i \rVert \leq C_{\theta} \leq \pi\ \forall\ \theta_i \in \Phi_i\ \text{and i} \in \small[N \small]
    \end{equation}
\end{assumption}
Assumptions \ref{ass:a1} and \ref{ass:a2} are generally considered for any standard online bounded-gradient-based Algorithm  so that adversaries can only exhibit limited manipulation of objectives.
\begin{assumption} \label{ass:a3}
The Lagrange multipliers $\{\lambda_{i,t}\}_{i \in \small[N \small],\ t \geq 1}$; associated with the power-balance constraint are bounded i.e., there exists a constant $C_{\lambda}$ s.t.\ $\forall\ i \in \small[N \small],\ t \geq 1$,
    \begin{equation}\label{eqn:e14}
    \lVert \lambda_{i,t} \rVert \leq C_{\lambda}
    \end{equation}
    Note $\lambda_{i,t}$ can represent the  Locational Marginal Price (LMP) at a node, and is a finite measure of node-specific local prices of delivering an additional unit of electricity within wholesale electricity markets \cite[Chap. 8]{raikar2019renewable}. Hence it is reasonable to assume it is bounded.
\end{assumption}

\begin{assumption} \label{ass:a4}
The weighted graph $\mathcal{G} = (\mathcal{V}, \mathcal{E}, W)$ satisfies: 
    \begin{enumerate}
        \item[a)] There exists a scalar $\eta \in (0,1)$ s.t. $W_{ij} \geq \eta$ if $j \in \mathcal{N}_i$, else $W_{ij} = 0$.
        \item[b)] Weight matrix  $W$ is doubly stochastic i.e. $\overset{N}{\underset{i=1}{\Sigma}} W_{ij} = 1\ \forall\ j \in \mathcal{V}$ and $\overset{N}{\underset{j=1}{\Sigma}} W_{ij} = 1\ \forall\ i \in \mathcal{V}$.
        \item[c)] 
        $\mathcal{G}(\mathcal{V}, \mathcal{E}, W)$ is undirected and connected.
    \end{enumerate}
\end{assumption}
 Assumption \ref{ass:a4} is needed to bound a certain sequence related to $\lambda_{i,t}$ which will be encountered later. 
\begin{assumption} \label{ass:a5}
Loads $\{p_i^l\ \forall\ i \in \small[N \small]\}$, present in the system, are considered to be time-invariant.
\end{assumption}

\subsection{General Relations}
In this module, we provide certain basic relations involving sequences generated by Algorithm \ref{algo:doopt}, which will play a crucial role in the subsequent analytical treatment of the sublinearly bounded regret and violation. 
\begin{lemma} \label{lem:l1}
Consider the sequences $\{p_{i,t}\}_{i \in \small[N \small],\ t \geq 1}$, $\{\theta_{i,t}\}_{i \in \small[N \small],\ t \geq 1}$ and $\{\lambda_{i,t}\}_{i \in \small[N \small],\ t \geq 1}$, as generated by Algorithm \ref{algo:doopt}. Define $N$-\textit{ary Cartesian product} over all sets $P_i$ and $\Phi_i$ as $P := \overset{N}{\underset{i=1}{\times}}P_i$, and $\Phi := \overset{N}{\underset{i=1}{\times}}\Phi_i$, respectively. Let $\mathbf{p_t^g} = \big[p_{1,t}, \hdots ,p_{N,t} \big]^T$,\ $\mathbf{\theta_t} = \big[\theta_{1,t}, \hdots ,\theta_{N,t} \big]^T$ and $\bar{\lambda}_t = \frac{1}{N} \overset{N}{\underset{i=1}{\Sigma}} \lambda_{i,t}$. Then, the following holds :
\begin{enumerate}[label=(\alph*)]
    \item For any $\mathbf{p^g} \in P$ and $\mathbf{\theta} \in \Phi$  over $t \geq 1$,
    \begin{equation}\label{eqn:e40}
        \begin{aligned}
            \mathcal{L}_t&([\mathbf{p_t^g}, \mathbf{\theta_t}], \bar{\lambda}_t) - \mathcal{L}_t([\mathbf{p^g}, \mathbf{\theta}], \bar{\lambda}_t)
            \\&\ 
            \leq \frac{1}{2\delta_t} \Big[\overset{N}{\underset{i=1}{\Sigma}} \big\lVert p_{i,t}^g - p_i^g \big\rVert^2 - \overset{N}{\underset{i=1}{\Sigma}} \lVert p_{i,t+1}^g - p_i^g \rVert^2 \Big] 
            \\&\
            + \frac{N \delta_t}{2} (L_f + C_{\lambda})^2 + 2C_P \overset{N}{\underset{i=1}{\Sigma}} \lVert \tilde{\lambda}_{i,t} - \bar{\lambda}_t \rVert
        \end{aligned}
    \end{equation}
    \item For any $\lambda \in \mathbb{R}$, and $t \geq 1$,
    \begin{equation*}
        \begin{aligned}
            \mathcal{L}_t&([\mathbf{p_t^g}, \mathbf{\theta_t}], \lambda) - \mathcal{L}_t([\mathbf{p_t^g}, \mathbf{\theta_t}], \bar{\lambda}_t)
            \\&\ 
            \leq \frac{1}{2\alpha_t} \Big[\overset{N}{\underset{i=1}{\Sigma}} \big\lVert \lambda_{i,t} - \lambda \big\rVert^2 - \overset{N}{\underset{i=1}{\Sigma}} \lVert \lambda_{i,t+1} - \lambda \rVert^2 \Big] +
        \end{aligned}
    \end{equation*}
    \begin{equation}\label{eqn:e41}
        \begin{aligned}
            &\
            (C_P + 2C_{\theta} \lVert B \rVert_{\infty}) \overset{N}{\underset{i=1}{\Sigma}} \lVert \tilde{\lambda}_{i,t} - \bar{\lambda}_t \rVert 
            + \frac{4N C_{\lambda}^2 \lVert B \rVert_{\infty} \beta_t}{\alpha_t} 
            \\&\
            + N(C_P + 2C_{\theta} \lVert B \rVert_{\infty})^2 \alpha_t + \frac{2N C_{\lambda}^2 \lVert B \rVert_{\infty}^2 \beta_t^2}{\alpha_t}
            \\&\ 
            + 2N C_{\lambda} \lVert B \rVert_{\infty} (C_P + 2C_{\theta} \lVert B \rVert_{\infty}) \beta_t 
        \end{aligned}
    \end{equation}
\end{enumerate}
\textit{Proof :} Refer to Section \ref{pflem1} of the Appendix for proof. \qed
\end{lemma}
\subsection{Main Results}
This section will be used to exhibit the sublinearity of the static regret and violation of constraints defined in \eqref{eqn:reg} and \eqref{eqn:ecviol}, for our devised algorithm.  
\begin{theorem}
Let Assumptions \ref{ass:a1} - \ref{ass:a5} be true. Consider the sequences $\{p_{i,t}\}_{i \in \small[N \small],\ t \geq 1}$, $\{\theta_{i,t}\}_{i \in \small[N \small],\ t \geq 1}$ and $\{\lambda_{i,t}\}_{i \in \small[N \small],\ t \geq 1}$, as generated by the Algorithm \eqref{algo:doopt}. Then, for the following choice of step size
\begin{equation*}
\alpha_t = \delta_t = \frac{1}{\sqrt{t}},\text{and } \beta_t = \frac{1}{t},
\end{equation*}
the static regret $\mathcal{R}_s(T)$ as defined in \eqref{eqn:reg} for any $T \geq 1$ satisfies,
\begin{equation*}
    \mathcal{R}_s(T) \leq K_1(N) + K_2(N) \sqrt{T} + K_3(N)\ ln(T) + K_4(N) \frac{1}{\sqrt{T}},
\end{equation*}
where the time-invariant parameters $K_1(N), K_2(N), K_3(N),$ and $K_4(N)$, has been defined in the \ref{nomen} of Appendix.
\end{theorem}
\textit{Proof :}
Invoking the results of \textit{Lemma 1}, by adding \eqref{eqn:e40} and \eqref{eqn:e41} and summing the expression over $t = 1, \hdots ,T$; we have for any $p^g \in P$ and $\lambda \in \mathbb{R}$,
\begin{equation*}
    \begin{aligned}
        &\ \overset{T}{\underset{t=1}{\Sigma}} \big[ \mathcal{L}_t([p_t^g, \theta_t], \lambda) - \mathcal{L}_t([p^g, \theta], \bar{\lambda}_t) \big]
        \\&\ 
        \leq \overset{T}{\underset{t=1}{\Sigma}} \frac{1}{2\delta_t} \Big[\overset{N}{\underset{i=1}{\Sigma}} \big\lVert p_{i,t}^g - p_i^g \big\rVert^2 - \overset{N}{\underset{i=1}{\Sigma}} \lVert p_{i,t+1}^g - p_i^g \rVert^2 \Big]\ + 
        \\&\
        \overset{T}{\underset{t=1}{\Sigma}} \frac{N \delta_t}{2} (L_f + C_{\lambda})^2 + 2C_P \overset{T}{\underset{t=1}{\Sigma}} \overset{N}{\underset{i=1}{\Sigma}} \lVert \tilde{\lambda}_{i,t} - \bar{\lambda}_t \rVert\ +
        \\&\
        \overset{T}{\underset{t=1}{\Sigma}} \frac{1}{2\alpha_t} \Big[\overset{N}{\underset{i=1}{\Sigma}} \big\lVert \lambda_{i,t} - \lambda \big\rVert^2 - \overset{N}{\underset{i=1}{\Sigma}} \lVert \lambda_{i,t+1} - \lambda \rVert^2 \Big]\ + 
        \\&\
        (C_P + 2C_{\theta} \lVert B \rVert_{\infty}) \overset{T}{\underset{t=1}{\Sigma}} \overset{N}{\underset{i=1}{\Sigma}} \lVert \tilde{\lambda}_{i,t} - \bar{\lambda}_t \rVert + \overset{T}{\underset{t=1}{\Sigma}} \frac{4N C_{\lambda}^2 \lVert B \rVert_{\infty} \beta_t}{\alpha_t}  
        \\&\
        + N(C_P + 2C_{\theta} \lVert B \rVert_{\infty})^2 \overset{T}{\underset{t=1}{\Sigma}} \alpha_t + \overset{T}{\underset{t=1}{\Sigma}} \frac{2N C_{\lambda}^2 \lVert B \rVert_{\infty}^2 \beta_t^2}{\alpha_t}\ + 
        \\&\
        2N C_{\lambda} \lVert B \rVert_{\infty} (C_P + 2C_{\theta} \lVert B \rVert_{\infty}) \overset{T}{\underset{t=1}{\Sigma}} \beta_t 
    \end{aligned}
\end{equation*}
\begin{equation*} 
    \begin{aligned}
        &\
        \overset{(a)}{\leq} \Big(\frac{2NC_P^2}{\delta_T} + \frac{2NC_{\lambda}^2}{\alpha_T} \Big) + \frac{N}{2} (L_f + C_{\lambda})^2 \small [2\sqrt{T} - 1\small]\ + 
        \\&\
        (3C_P + 2C_{\theta} \lVert B \rVert_\infty) \overset{T}{\underset{t=1}{\Sigma}} \overset{N}{\underset{i=1}{\Sigma}} \lVert \tilde{\lambda}_{i,t} - \bar{\lambda}_t \rVert + 4N C_{\lambda}^2 \lVert B \rVert_{\infty}  
        \\&\
        [2\sqrt{T} - 1\small] + N(C_P + 2C_{\theta} \lVert B \rVert_{\infty})^2 [2\sqrt{T} - 1\small] + 2N C_{\lambda}^2 \lVert B \rVert_{\infty}^2  
        \\&\
        [3 - \frac{2}{\sqrt{T}}\small] + 2N C_{\lambda} \lVert B \rVert_{\infty} (C_P + 2C_{\theta} \lVert B \rVert_{\infty})\ \small[1 + ln(T)\small]
        \\&\ 
        \overset{(b)}{\leq} \Big\{2N(C_P^2 + C_{\lambda}^2) + N(L_f + C_{\lambda})^2\ +
        \\&\
        8N C_{\lambda}^2 \lVert B \rVert_{\infty} + 2N(C_P + 2C_{\theta} \lVert B \rVert_{\infty})^2\ + 
        \\&\
        2(3C_P + 2C_{\theta} \lVert B \rVert_\infty)(C_P + 2C_{\theta} \lVert B \rVert_{\infty})(2N + \kappa N^2) \Big\}\sqrt{T}\ - 
        \\&\
        \big\{4N C_{\lambda}^2 \lVert B \rVert_{\infty}^2\big\} \frac{1}{\sqrt{T}} + \Big\{2N C_{\lambda} \lVert B \rVert_{\infty} (C_P + 2C_{\theta} \lVert B \rVert_{\infty})\ + 
    \end{aligned}
\end{equation*}
\begin{equation}\label{eqn:e42}
    \begin{aligned}    
        &\
        (3C_P + 2C_{\theta} \lVert B \rVert_\infty)(4NC_{\lambda} \lVert B \rVert_\infty + 2\kappa N^2 C_{\lambda} \lVert B \rVert_{\infty} \overset{T-2}{\underset{l=1}{\Sigma}} \omega^l) \Big\}
        \\&\
        ln(T) - \Big\{\frac{N}{2} (L_f + C_{\lambda})^2 +  N(C_P + 2C_{\theta} \lVert B \rVert_{\infty})^2\ - 
        \\&\
        6N C_{\lambda}^2 \lVert B \rVert_{\infty}^2 - 2N C_{\lambda} \lVert B \rVert_{\infty} (C_P + 2C_{\theta} \lVert B \rVert_{\infty})\ + 
        \\&\
        4N C_{\lambda}^2 \lVert B \rVert_{\infty} - \big(3C_P\ + 2C_{\theta} \lVert B \rVert_\infty\big) \Big[2NC_{\lambda} + \frac{\kappa N^2 \omega C_{\lambda}}{1-\omega}
        \\&\
        - 2N(C_P + 2C_{\theta} \lVert B \rVert_{\infty})\ + 4NC_{\lambda} \lVert B \rVert_\infty\ -
        \\&\
        \kappa N^2 (C_P + 2C_{\theta} \lVert B \rVert_{\infty}) \overset{T-2}{\underset{l=1}{\Sigma}} \omega^l + 2\kappa N^2 C_{\lambda} \lVert B \rVert_{\infty} \overset{T-2}{\underset{l=1}{\Sigma}} \omega^l \Big]\Big\}
        \\&\
        = K_1(N) + K_2(N) \sqrt{T} + K_3(N)\ ln(T) + K_4(N) \frac{1}{\sqrt{T}}
    \end{aligned}
\end{equation}
\textit{where, (a) follows from \eqref{eqn:e24}, \eqref{eqn:e25}, \eqref{eqn:e32}, \eqref{eqn:e33} \& Corollary \ref{cor:c7} and, (b) follows from Corollary \ref{cor:c9}.} 

Since, the inequality \eqref{eqn:e42} holds true for any $p^g \in P$ and, $\lambda \in \mathbb{R}$, let $p^g \triangleq p^{g*} \in P$ and $\lambda \triangleq \textbf{0}$. From this, we obtain,
\begin{equation}\label{eqn:e43}
    \begin{aligned}
        \mathcal{R}_s(T) &\ = \overset{T}{\underset{t=1}{\Sigma}} \big[ \mathcal{L}_t([p_t^g, \theta_t], \textbf{0}) - \mathcal{L}_t([p^{g*}, \theta^{*}], \bar{\lambda}_t) \big]
        \\&\
        = \overset{T}{\underset{t=1}{\Sigma}} \big[\overset{N}{\underset{i=1}{\Sigma}}f_{i,t}(p_{i,t}^g) - \overset{N}{\underset{i=1}{\Sigma}}f_{i,t}(p_i^{g*}) - \overset{N}{\underset{i=1}{\Sigma}}\bar{\lambda}_t h_i(p_i^{g*}, \theta^{*}) \big]
        \\&\
        \overset{(a)}{=} \overset{T}{\underset{t=1}{\Sigma}} \overset{N}{\underset{i=1}{\Sigma}} \big[f_{i,t}(p_{i,t}^g) - f_{i,t}(p_i^{g*}) \big]
        \\&\
        \overset{(b)}{\leq} K_1(N) + K_2(N) \sqrt{T} + K_3(N)\ ln(T) + K_4(N) \frac{1}{\sqrt{T}}
    \end{aligned} 
\end{equation}
\textit{where, (a) borrows from \eqref{eqn:e3} and primal feasibility of KKT conditions and, (b) follows from \eqref{eqn:e42}.} \qed

\begin{theorem}
Let Assumptions \ref{ass:a1} - \ref{ass:a5} be true. Consider the sequences $\{p_{i,t}\}_{i \in \small[N \small],\ t \geq 1}$, $\{\theta_{i,t}\}_{i \in \small[N \small],\ t \geq 1}$ and $\{\lambda_{i,t}\}_{i \in \small[N \small],\ t \geq 1}$, as generated by the Algorithm \eqref{algo:doopt}. Then, for the following choice of step size
\begin{equation*}
    \gamma_t = \frac{1}{\sqrt{t}},
\end{equation*}
the constraint violation $\mathcal{R}_s^{ec}(T)$ as defined in \eqref{eqn:ecviol} for any $T \geq 1$ satisfies,
\begin{equation*}
    \mathcal{R}_s^{ec}(T) \leq M_1(N)\sqrt{T},
\end{equation*}
where the time-invariant parameter $M_1(N)$, has been defined in the \ref{nomen} of Appendix.
\end{theorem}

\textit{Proof :}
Using \eqref{eqn:e8},
\begin{equation}
    \begin{aligned}
        \theta_{i,t+1} &\ = \theta_{i,t} - \gamma_t \nabla_{\lambda_{i,t}}\mathcal{L}_t
        \\&\
        = \theta_{i,t} - \gamma_t \big[p^g_{i,t} - p^l_{i,t} - \theta_{i,t}\underset{j \in \mathcal{N}_i}{\Sigma}B_{ij} + \underset{j \in \mathcal{N}_i}{\Sigma}B_{ij} \theta_{j,t} \big]
    \end{aligned}
\end{equation}
\big[Summing over all agents\big]
\begin{equation*}
    \begin{aligned}
        \overset{N}{\underset{i=1}{\Sigma}} \theta_{i,t+1} =&\ \overset{N}{\underset{i=1}{\Sigma}} \theta_{i,t} - \overset{N}{\underset{i=1}{\Sigma}} \gamma_t \big[p^g_{i,t} - p^l_{i,t} - \theta_{i,t} \underset{j \in \mathcal{N}_i}{\Sigma}B_{ij}\ + 
        \\&\
        \underset{j \in \mathcal{N}_i}{\Sigma}B_{ij} \theta_{j,t}\big]
    \end{aligned}
\end{equation*}
\begin{equation}
    \begin{aligned}
        or,\ \overset{N}{\underset{i=1}{\Sigma}} \big[&p^g_{i,t} - p^l_{i,t} - \theta_{i,t}\underset{j \in \mathcal{N}_i}{\Sigma}B_{ij} + \underset{j \in \mathcal{N}_i}{\Sigma}B_{ij} \theta_{j,t}\big]\ = 
        \\&\
        \frac{\overset{N}{\underset{i=1}{\Sigma}} \theta_{i,t+1} - \overset{N}{\underset{i=1}{\Sigma}} \theta_{i,t}}{\gamma_t} \leq \frac{\overset{N}{\underset{i=1}{\Sigma}} \theta_{i,t+1} - \overset{N}{\underset{i=1}{\Sigma}} \theta_{i,t}}{\gamma_T}
    \end{aligned}
\end{equation}
\big[Summing over all time\big]
\begin{equation}
    \begin{aligned}
        &\ \overset{T}{\underset{t=1}{\Sigma}} \overset{N}{\underset{i=1}{\Sigma}} \big[p^g_{i,t} - p^l_{i,t} - \theta_{i,t}\underset{j \in \mathcal{N}_i}{\Sigma}B_{ij} + \underset{j \in \mathcal{N}_i}{\Sigma}B_{ij} \theta_{j,t}\big]\ \leq 
        \\&\
        \frac{\overset{T}{\underset{t=1}{\Sigma}} \overset{N}{\underset{i=1}{\Sigma}} \theta_{i,t+1} - \overset{T}{\underset{t=1}{\Sigma}} \overset{N}{\underset{i=1}{\Sigma}} \theta_{i,t}}{\gamma_T} = \frac{\overset{N}{\underset{i=1}{\Sigma}} \theta_{i,T+1}}{\gamma_T}
    \end{aligned}
\end{equation}
Hence, we get,
\begin{equation}
    \begin{aligned}
        \mathcal{R}_s^{ec}(T) &\ = \Big\lVert \overset{T}{\underset{t=1}{\Sigma}} \overset{N}{\underset{i=1}{\Sigma}} p^g_{i,t} - p^l_{i,t} - \theta_{i,t}\underset{j \in \mathcal{N}_i}{\Sigma}B_{ij} + \underset{j \in \mathcal{N}_i}{\Sigma}B_{ij} \theta_{j,t} \Big\rVert 
        \\&\
         \leq \frac{NC_{\theta}}{\gamma_T} \overset{(a)}{=} NC_{\theta}\sqrt{T} = M_1(N)\sqrt{T}
    \end{aligned}
\end{equation}
\textit{where, (a) implies choosing $\gamma_t = \frac{1}{\sqrt{t}}$.} \qed

\section{Numerical Example} \label{numex}
In this section, we use the standard IEEE-14 bus test case \cite{aeps1962data} to demonstrate the convergence of Algorithm \ref{algo:doopt}.
\begin{figure} 
    \centerline{\includegraphics[scale=0.6]{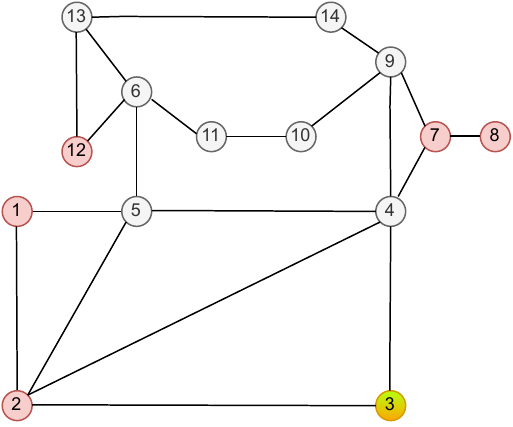}}
   \caption{IEEE 14-bus system. The nodes of the network represent the buses with red ones being generators, yellow being slack and a few being loads.}
   \label{fig:grp}
\end{figure}

There are $N=14$ agents interacting throughout the system over an undirected network as shown in Fig. \eqref{fig:grp}. The power transmission network topology and communication topology are considered to be the same in this case. Out of the 14 buses, generators are located on five of them viz. Nodes 1, 2, 7, 8 and 12. Node 3 has been designated to be the slack bus. The generator limits have been provided in Table \ref{tab:table}. Some buses are profiled with constant-power load viz. Nodes 3, 4, 6, 8, 10, 12, and 14 are having 95, 125, 100, 106, 115, 90, and 80 MW as loads, respectively.  
\begin{table}[ht!] 
  \begin{center}
    \begin{tabular}{c|c|c}
      \hline 
      Location & $P^g_{i,min} (MW)$ & $P^g_{i,min} (MW)$\\
      \hline
      Node 1 & 0 & 280\\
      Node 2 & 0 & 110\\
      Node 7 & 0 & 70\\
      Node 8 & 0 & 160\\
      Node 12 & 0 & 130\\
      \hline 
    \end{tabular}
    \vspace{0.3cm}
    \caption{Generator Parameters}
    \label{tab:table}
  \end{center}
\end{table}

A time horizon of $T=2000$ instances has been considered and for the objective function $f_{i,t}(p^g_{i,t}) = a^t_i(p^g_{i,t})^2 + b^t_i(p^g_{i,t}) + c^t_i$, we have generated synthetic data by drawing from the uniform distribution, $a^t_i \in \mathcal{U}_{\small[0.001, 0.08\small]}$, $b^t_i \in \mathcal{U}_{\small[1, 5\small]}$ and $c^t_i$ has been set to all zeros.
\begin{figure} 
  \centerline{\includegraphics[scale=0.3]{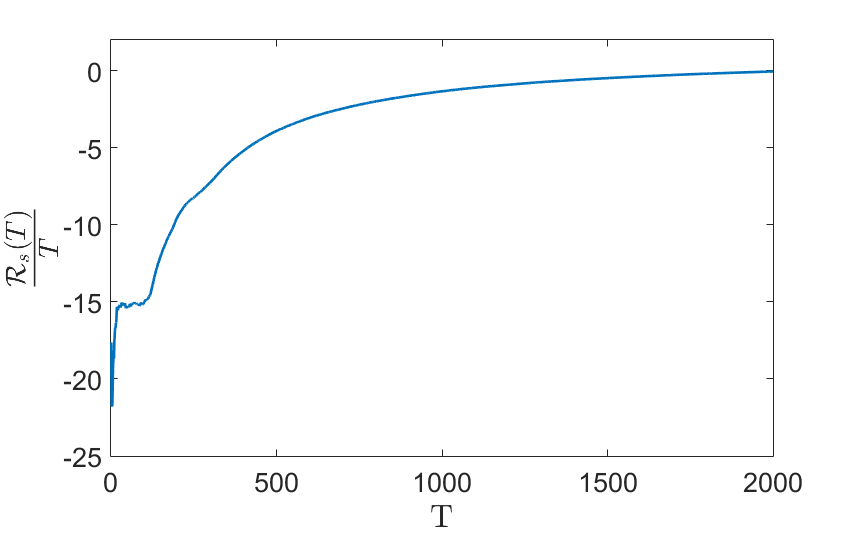}}
   \caption{Dynamics of the Average Regret against time-horizon for $N = 14$.}
   \label{fig:reg}
\end{figure}
\begin{figure} 
  \centerline{\includegraphics[scale=0.3]{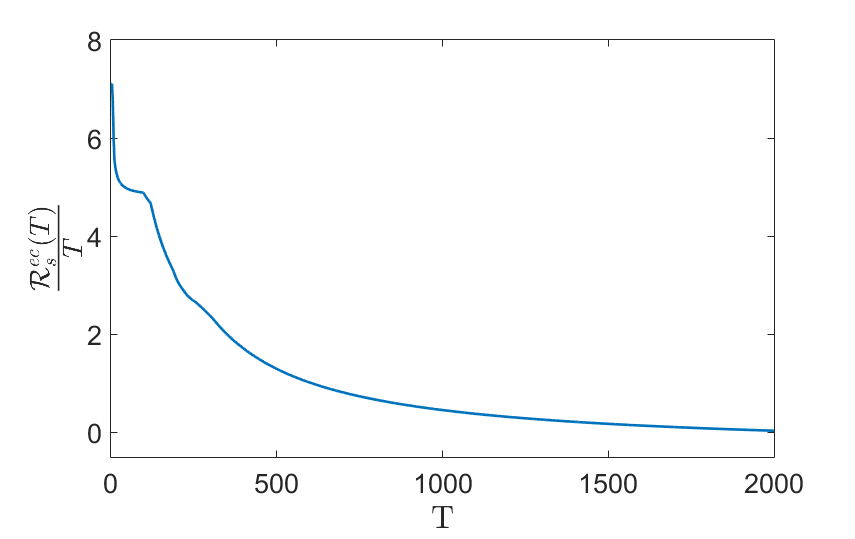}}
   \caption{Dynamics of the Average Violation against time-horizon for $N = 14$.}
   \label{fig:viol}
\end{figure}
\begin{remark}
Depending on the conditioning of initial values and learning rate parameters used in the algorithm, the regret can assume negative values \cite{li2020distributed} however, it doesn't necessarily imply that the online solution is "better" than the offline one initially, just that the algorithm \text{greedily} starts to explore the decision space with past information. Irrespectively, the primary objective of achieving \textit{sublinear} regret over the horizon remains the ultimate goal and as such the online solution should "catch up" with the offline one in the \textit{long-term}.   
\end{remark}
Fig. \eqref{fig:reg} and Fig. \eqref{fig:viol} illustrates the trajectories of the Average \textit{Static} Regret $\frac{\mathcal{R}_s(T)}{T}$ and Average Constraint Violation $\frac{\mathcal{R}^{ec}_s(T)}{T}$ for $N=14$ agents over $T=2000$ time horizon instances. It can be observed that Algorithm \ref{algo:doopt} ensures that both the performance metrics approach zero, as the horizon increases. In a sense, it can be inferred from this example that the Algorithm \ref{algo:doopt} is learning the cost parameter probability distribution of the objective functions, and as the time horizon increases so do the prediction capabilities.
\begin{figure} 
  \centerline{\includegraphics[scale=0.3]{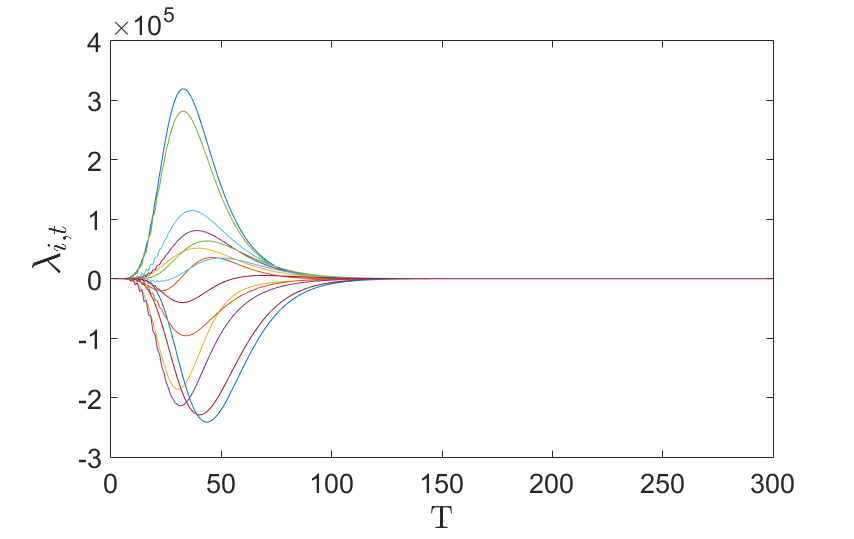}}
   \caption{Trajectories of the Locational Marginal Prices for each agent; $\lambda_{i,t}$.}
   \label{fig:mu}
\end{figure}
The evolution of the Localized Marginal Prices (LMP) $\lambda_{i,t}$ for each agent $i \in \small[N\small]$ has been depicted in Fig. \eqref{fig:mu}.
\begin{remark}
The primary objective of the online OPF algorithm is to learn the uncertainties inherent in the optimization problem so as to make better and more informed decisions regarding optimized generation allocation under inconclusiveness, and as such is usually operated initially in the background \cite{chadoulos2021learning}, and deployed in real-time once training finishes.   
\end{remark}

\section{Conclusion} \label{concl}
    In this paper, the distributed Optimal Power Flow (OPF) problem in an online learning-based environment, has been investigated. Each agent in this power network is only privy to its local objective function and constraints, which are dynamic in nature and revealed to the agents only after a decision has been made at each stage. 
    To this end, a distributed online algorithm based on the primal-dual approach has been proposed to solve this problem. Theoretical guarantees on the convergence of the algorithm have been provided through the sublinear bounding of static regret and constraint violation. Finally, using numerical illustrations, the performance of the algorithm has been demonstrated. For future work, the algorithm in this paper can be extended to account for the line congestion and stability constraints.

\section*{Appendix} \label{appdx}
\subsection{Nomenclature} \label{nomen}
\begin{equation*} 
    \begin{aligned}
        K_1(N) =& -\frac{N}{2} (L_f + C_{\lambda})^2 -  N(C_P + 2C_{\theta} \lVert B \rVert_{\infty})^2\ +  
        \\&\
        6N C_{\lambda}^2 \lVert B \rVert_{\infty}^2 + 2N C_{\lambda} \lVert B \rVert_{\infty} (C_P + 2C_{\theta} \lVert B \rVert_{\infty})\ - 
        \\&\
        4N C_{\lambda}^2 \lVert B \rVert_{\infty} + (3C_P + 2C_{\theta} \lVert B \rVert_\infty) \Big[2NC_{\lambda} - \frac{\kappa N^2 \omega C_{\lambda}}{1-\omega} 
        \\&\
        + 2N(C_P + 2C_{\theta} \lVert B \rVert_{\infty}) - 2\kappa N^2 C_{\lambda} \lVert B \rVert_{\infty} \overset{T-2}{\underset{l=1}{\Sigma}} \omega^l\ +
        \\&\
        \kappa N^2 (C_P + 2C_{\theta} \lVert B \rVert_{\infty}) \overset{T-2}{\underset{l=1}{\Sigma}} \omega^l \Big] - 4NC_{\lambda} \lVert B \rVert_\infty\\
        K_2(N) =&\ 2N(C_P^2 + C_{\lambda}^2) + N(L_f + C_{\lambda})^2 + 8N C_{\lambda}^2 \lVert B \rVert_{\infty}\ + 
        \\&\
        2N(C_P + 2C_{\theta} \lVert B \rVert_{\infty})^2\ + 
        \\&\
        2(3C_P + 2C_{\theta} \lVert B \rVert_\infty)(C_P + 2C_{\theta} \lVert B \rVert_{\infty})(2N + \kappa N^2)\\
        K_3(N) = &\ 2N C_{\lambda} \lVert B \rVert_{\infty} (C_P + 2C_{\theta} \lVert B \rVert_{\infty}) + (3C_P + 2C_{\theta} \lVert B \rVert_\infty)
        \\&\
        (4NC_{\lambda} \lVert B \rVert_\infty + 2\kappa N^2 C_{\lambda} \lVert B \rVert_{\infty} \overset{T-2}{\underset{l=1}{\Sigma}} \omega^l)\\
        K_4(N) =&\ -4N C_{\lambda}^2 \lVert B \rVert_{\infty}^2 \\
        M_1(N) =&\ NC_{\theta}
    \end{aligned} 
\end{equation*}
 
\subsection{Proof of Lemma 1}\label{pflem1}
\begin{enumerate}[label=(\alph*)]
    \item Let Assumption \ref{ass:a5} hold. Now, for any $p^g \in P$, $\theta \in \Phi$ and $t \geq 1$,
    \begin{equation*}
        \begin{aligned}
            &\ \mathcal{L}_t([p_t^g, \theta_t], \bar{\lambda}_t) - \mathcal{L}_t([p^g, \theta], \bar{\lambda}_t)
            \\&\ 
            = \overset{N}{\underset{i=1}{\Sigma}} f_{i,t}(p^g_{i,t}) + \overset{N}{\underset{i=1}{\Sigma}} \bar{\lambda}_t (p^g_{i,t} - p^l_{i,t} - \theta_{i,t}\underset{j \in \mathcal{N}_i}{\Sigma}B_{ij}\ + 
            \\&\
            \underset{j \in \mathcal{N}_i}{\Sigma}B_{ij} \theta_{j,t}) - \overset{N}{\underset{i=1}{\Sigma}} f_{i,t}(p^g_i) - \overset{N}{\underset{i=1}{\Sigma}} \bar{\lambda}_t (p^g_i - p^l_i\ - 
        \end{aligned}
    \end{equation*}
    \begin{equation}\label{eqn:e80}
        \begin{aligned}
            &\
            \theta_i \underset{j \in \mathcal{N}_i}{\Sigma}B_{ij} + \underset{j \in \mathcal{N}_i}{\Sigma}B_{ij} \theta_j) 
            \\&\ 
            \overset{(a)}{=} \overset{N}{\underset{i=1}{\Sigma}} \small[f_{i,t}(p^g_{i,t}) - f_{i,t}(p^g_i)\small] + \bar{\lambda}_t \Big\{\overset{N}{\underset{i=1}{\Sigma}} \small[p^g_{i,t} - p^g_i \small]\ +
            \\&\
            \overset{N}{\underset{i=1}{\Sigma}} \small[\theta_i - \theta_{i,t}\small] \underset{j \in \mathcal{N}_i}{\Sigma} B_{ij} + \overset{N}{\underset{i=1}{\Sigma}} \underset{j \in \mathcal{N}_i}{\Sigma} B_{ij} \small[\theta_{j,t} - \theta_j\small] \Big\}\ +  
            \\&\ 
            \frac{1}{2\delta_t} \Big(\overset{N}{\underset{i=1}{\Sigma}} \lVert p_{i,t+1}^g - p_i^g \rVert^2 - \overset{N}{\underset{i=1}{\Sigma}} \lVert p_{i,t+1}^g - p_i^g \rVert^2\Big) 
            \\&\
            \overset{(b)}{=} \overset{N}{\underset{i=1}{\Sigma}} \small[f_{i,t}(p^g_{i,t}) - f_{i,t}(p^g_i)\small] + \bar{\lambda}_t \Big\{\overset{N}{\underset{i=1}{\Sigma}} \small[p^g_{i,t} - p^g_i \small]\ + 
            \\&\
            \overset{N}{\underset{i=1}{\Sigma}} \small[\theta_i - \theta_{i,t}\small] \underset{j \in \mathcal{N}_i}{\Sigma} B_{ij} + \overset{N}{\underset{i=1}{\Sigma}} \underset{j \in \mathcal{N}_i}{\Sigma} B_{ij} \small[\theta_{j,t} - \theta_j\small] \Big\} - \frac{1}{2\delta_t}
            \\&\
            \Big(\overset{N}{\underset{i=1}{\Sigma}} \lVert p_{i,t+1}^g - p_i^g \rVert^2 + \overset{N}{\underset{i=1}{\Sigma}} \big\lVert \mathbb{P}_{P_i} \big[p_{i,t}^g - \delta_t \nabla_{p^g_{i,t}}\mathcal{L}_t\big] - p_i^g \big\rVert^2 \Big) 
            \\&\
            \overset{(c)}{\leq} \overset{N}{\underset{i=1}{\Sigma}} \small[f_{i,t}(p^g_{i,t}) - f_{i,t}(p^g_i)\small] + \bar{\lambda}_t \Big\{\overset{N}{\underset{i=1}{\Sigma}} \small[p^g_{i,t} - p^g_i \small]\ + 
            \\&\
            \overset{N}{\underset{i=1}{\Sigma}} \small[\theta_i - \theta_{i,t}\small] \underset{j \in \mathcal{N}_i}{\Sigma} B_{ij} + \overset{N}{\underset{i=1}{\Sigma}} \underset{j \in \mathcal{N}_i}{\Sigma} B_{ij} \small[\theta_{j,t} - \theta_j\small] \Big\}\ -
            \\&\
            \frac{1}{2\delta_t} \overset{N}{\underset{i=1}{\Sigma}} \lVert p_{i,t+1}^g - p_i^g \rVert^2 + \frac{1}{2\delta_t} \overset{N}{\underset{i=1}{\Sigma}} \big\lVert p_{i,t}^g - \delta_t \nabla_{p^g_{i,t}}\mathcal{L}_t - p_i^g \big\rVert^2
            \\&\
            \overset{(d)}{\leq} \frac{1}{2\delta_t} \Big[\overset{N}{\underset{i=1}{\Sigma}} \big\lVert p_{i,t}^g - p_i^g \big\rVert^2 - \overset{N}{\underset{i=1}{\Sigma}} \lVert p_{i,t+1}^g - p_i^g \rVert^2 \Big]\ +
            \\&\
            \frac{N \delta_t}{2} (L_f + C_{\lambda})^2 - \overset{N}{\underset{i=1}{\Sigma}} \tilde{\lambda}_{i,t} \small[p^g_{i,t} - p^g_i \small]
            + \bar{\lambda}_t \Big\{\overset{N}{\underset{i=1}{\Sigma}} \small[p^g_{i,t} - p^g_i \small]
            \\&\
            + \overset{N}{\underset{i=1}{\Sigma}} \small[\theta_i - \theta_{i,t}\small] \underset{j \in \mathcal{N}_i}{\Sigma} B_{ij} + \overset{N}{\underset{i=1}{\Sigma}} \underset{j \in \mathcal{N}_i}{\Sigma} B_{ij} \small[\theta_{j,t} - \theta_j\small] \Big\}
            \\&\
            \overset{(e)}{=} \frac{1}{2\delta_t} \Big[\overset{N}{\underset{i=1}{\Sigma}} \big\lVert p_{i,t}^g - p_i^g \big\rVert^2 - \overset{N}{\underset{i=1}{\Sigma}} \lVert p_{i,t+1}^g - p_i^g \rVert^2 \Big]\ + 
            \\&\
            \frac{N \delta_t}{2} (L_f + C_{\lambda})^2 - \overset{N}{\underset{i=1}{\Sigma}} \tilde{\lambda}_{i,t} \small[p^g_{i,t} - p^g_i \small]
            + \bar{\lambda}_t \overset{N}{\underset{i=1}{\Sigma}} \small[p^g_{i,t} - p^g_i \small]
            \\&\
            + \Big(\overset{N}{\underset{i=1}{\Sigma}} \bar{\lambda}_t \small[p^g_{i,t} - p^g_i \small] - \overset{N}{\underset{i=1}{\Sigma}} \bar{\lambda}_t \small[p^g_{i,t} - p^g_i \small]\Big)
            \\&\
            \overset{(f)}{=} \frac{1}{2\delta_t} \Big[\overset{N}{\underset{i=1}{\Sigma}} \big\lVert p_{i,t}^g - p_i^g \big\rVert^2 - \overset{N}{\underset{i=1}{\Sigma}} \lVert p_{i,t+1}^g - p_i^g \rVert^2 \Big]\ + 
            \\&\
            \frac{N \delta_t}{2} (L_f + C_{\lambda})^2 - \overset{N}{\underset{i=1}{\Sigma}} (\tilde{\lambda}_{i,t} - \bar{\lambda}_t) \small[p^g_{i,t} - p^g_i \small]
            \\&\
            \overset{(g)}{\leq} \frac{1}{2\delta_t} \Big[\overset{N}{\underset{i=1}{\Sigma}} \big\lVert p_{i,t}^g - p_i^g \big\rVert^2 - \overset{N}{\underset{i=1}{\Sigma}} \lVert p_{i,t+1}^g - p_i^g \rVert^2 \Big]\ + 
            \\&\
            \frac{N \delta_t}{2} (L_f + C_{\lambda})^2 + 2C_P \overset{N}{\underset{i=1}{\Sigma}} \lVert \tilde{\lambda}_{i,t} - \bar{\lambda}_t \rVert
        \end{aligned}
    \end{equation}
    \textit{where (a) implies introducing the last two self-cancelling terms, (b) follows from \eqref{eqn:e7} in the last term, and Assumption \ref{ass:a5} causes the third-last term to cancel out, (c) follows from the non-expansiveness property of euclidean projection operator $\mathbb{P}_{P_i} \big[. \big]$, (d) causes the first and second-last terms to cancel out in the previous expression and draws from \eqref{eqn:e23}, (e) involves \eqref{eqn:e26} and introducing the last two self-cancelling terms, (f) involves cancelling out last and third-last term in the previous expression, and (g) involves Assumption \ref{ass:a2}.}
    
    \item For any $\lambda \in \mathbb{R}$, and $t \geq 1$,
    \begin{equation*}
        \begin{aligned}
            &\ \mathcal{L}_t([\mathbf{p_t^g}, \mathbf{\theta_t}], \lambda) - \mathcal{L}_t([\mathbf{p_t^g}, \mathbf{\theta_t}], \bar{\lambda}_t)
            \\&\ 
            = \overset{N}{\underset{i=1}{\Sigma}} f_{i,t}(p^g_{i,t})\ +\ \overset{N}{\underset{i=1}{\Sigma}} \lambda (p^g_{i,t} - p^l_{i,t} - \theta_{i,t}\underset{j \in \mathcal{N}_i}{\Sigma}B_{ij}\ + 
        \end{aligned}
    \end{equation*}
    \begin{equation}\label{eqn:e81}
        \begin{aligned}
            &\
            \underset{j \in \mathcal{N}_i}{\Sigma}B_{ij} \theta_{j,t})
            - \overset{N}{\underset{i=1}{\Sigma}} f_{i,t}(p^g_{i,t})
            - \overset{N}{\underset{i=1}{\Sigma}} \bar{\lambda}_t \big(p^g_{i,t} - p^l_{i,t}\ - 
            \\&\
            \theta_{i,t}\underset{j \in \mathcal{N}_i}{\Sigma}B_{ij} + \underset{j \in \mathcal{N}_i}{\Sigma}B_{ij} \theta_{j,t}\big) 
            \\&\ 
            \overset{(a)}{=} \overset{N}{\underset{i=1}{\Sigma}} \small[\lambda - \bar{\lambda}_t\small] (p^g_{i,t} - p^l_{i,t} - \theta_{i,t}\underset{j \in \mathcal{N}_i}{\Sigma}B_{ij} + \underset{j \in \mathcal{N}_i}{\Sigma}B_{ij} \theta_{j,t})
            \\&\
            + \frac{1}{2 \alpha_t} \Big(\overset{N}{\underset{i=1}{\Sigma}} \lVert \lambda_{i,t+1} - \lambda \rVert^2 - \overset{N}{\underset{i=1}{\Sigma}} \lVert \lambda_{i,t+1} - \lambda \rVert^2\Big)
            \\&\
            \overset{(b)}{=} \overset{N}{\underset{i=1}{\Sigma}} \small[\lambda - \bar{\lambda}_t\small] (p^g_{i,t} - p^l_{i,t} - \theta_{i,t}\underset{j \in \mathcal{N}_i}{\Sigma}B_{ij} + \underset{j \in \mathcal{N}_i}{\Sigma}B_{ij} \theta_{j,t}) 
            \\&\
            - \frac{1}{2 \alpha_t} \overset{N}{\underset{i=1}{\Sigma}} \lVert \lambda_{i,t+1} - \lambda \rVert^2
            + \frac{1}{2 \alpha_t} \overset{N}{\underset{i=1}{\Sigma}} \lVert \tilde{\lambda}_{i,t} - \beta_t \nabla_{\theta_{i,t}}\mathcal{L}_t\ +
            \\&\
            \alpha_t \nabla_{\lambda_{i,t}}\mathcal{L}_t - \lambda \rVert^2
            \\&\
            \overset{(c)}{=} \overset{N}{\underset{i=1}{\Sigma}} \nabla_{\lambda_{i,t}}\mathcal{L}_t \small[\lambda - \bar{\lambda}_t\small] - \frac{1}{2 \alpha_t} \overset{N}{\underset{i=1}{\Sigma}} \lVert \lambda_{i,t+1} - \lambda \rVert^2\ + 
            \\&\
            \frac{1}{2 \alpha_t} \overset{N}{\underset{i=1}{\Sigma}} \lVert \tilde{\lambda}_{i,t} - \lambda \rVert^2 + \frac{1}{2 \alpha_t} \overset{N}{\underset{i=1}{\Sigma}} \lVert \alpha_t \nabla_{\lambda_{i,t}}\mathcal{L}_t - \beta_t \nabla_{\theta_{i,t}}\mathcal{L}_t \rVert^2
            \\&\
            + \frac{2}{2 \alpha_t} \overset{N}{\underset{i=1}{\Sigma}} \small[ \alpha_t \nabla_{\lambda_{i,t}}\mathcal{L}_t - \beta_t \nabla_{\theta_{i,t}}\mathcal{L}_t \small] (\tilde{\lambda}_{i,t} - \lambda)
            \\&\
            \overset{(d)}{=} \frac{1}{2\alpha_t} \Big[\overset{N}{\underset{i=1}{\Sigma}} \big\lVert \lambda_{i,t} - \lambda \big\rVert^2 - \overset{N}{\underset{i=1}{\Sigma}} \lVert \lambda_{i,t+1} - \lambda \rVert^2 \Big]\ + 
            \\&\
            \overset{N}{\underset{i=1}{\Sigma}} \nabla_{\lambda_{i,t}}\mathcal{L}_t \small[\tilde{\lambda}_{i,t} - \bar{\lambda}_t\small] - \frac{2\beta_t}{2\alpha_t} \overset{N}{\underset{i=1}{\Sigma}} \nabla_{\theta_{i,t}}\mathcal{L}_t \small[\tilde{\lambda}_{i,t} - \lambda \small]\ +
            \\&\ 
            \overset{N}{\underset{i=1}{\Sigma}} \frac{\alpha_t^2}{2 \alpha_t} \lVert \nabla_{\lambda_{i,t}}\mathcal{L}_t \rVert^2 + \overset{N}{\underset{i=1}{\Sigma}} \frac{\beta_t^2}{2 \alpha_t}  \lVert \nabla_{\theta_{i,t}}\mathcal{L}_t \rVert^2\ - 
            \\&\
            \overset{N}{\underset{i=1}{\Sigma}}
            \frac{2 \alpha_t \beta_t}{2 \alpha_t}  (\nabla_{\lambda_{i,t}}\mathcal{L}_t)(\nabla_{\theta_{i,t}}\mathcal{L}_t)
            \\&\
            \overset{(e)}{\leq} \frac{1}{2\alpha_t} \Big[\overset{N}{\underset{i=1}{\Sigma}} \big\lVert \lambda_{i,t} - \lambda \big\rVert^2 - \overset{N}{\underset{i=1}{\Sigma}} \lVert \lambda_{i,t+1} - \lambda \rVert^2 \Big] + 
            \\&\
            (C_P + 2C_{\theta} \lVert B \rVert_{\infty}) \overset{N}{\underset{i=1}{\Sigma}} \lVert \tilde{\lambda}_{i,t} - \bar{\lambda}_t \rVert + \frac{4N C_{\lambda}^2 \lVert B \rVert_{\infty} \beta_t}{\alpha_t}\ +
            \\&\
            N(C_P + 2C_{\theta} \lVert B \rVert_{\infty})^2 \alpha_t + \frac{2N C_{\lambda}^2 \lVert B \rVert_{\infty}^2 \beta_t^2}{\alpha_t}\ +
            \\&\ 
            2N C_{\lambda} \lVert B \rVert_{\infty} (C_P + 2C_{\theta} \lVert B \rVert_{\infty}) \beta_t
        \end{aligned}
    \end{equation}
    \textit{where (a) involves cancelling out first and third terms in the previous expression and introducing the last two self-cancelling terms, (b) follows from \eqref{eqn:e10}, (c) involves a non-stationary version of \eqref{eqn:e3} and expansion of the last term, (d) follows from cancelling terms and expanding the second-last term in the previous expression and (e) involves \eqref{eqn:e30} and \eqref{eqn:e31}.}
\end{enumerate}

\subsection{Corollaries and Lemmas}\label{corlem}

\begin{corollary} \label{cor:c1} 
Compactness of $P_i$ and continuity of $\{f_{i,t}\}_{i \in \small[N \small],\ t \geq 1}$ from Assumptions \ref{ass:a1} and \ref{ass:a2} yields a result that there exists a constant $C_f > 0$ s.t. $\forall\ i \in \small[N \small]$,
    \begin{equation}\label{eqn:e15}
    \lVert f_{i,t}(p^g_{i,t})  \rVert \leq C_f,\ \forall\ p^g_{i,t} \in P_i
    \end{equation}
\end{corollary}

\begin{corollary} \label{cor:c2}
Using Assumption \ref{ass:a1}, gradient boundedness and lipschitz continuity of $\{f_{i,t}\}_{i \in \small[N \small],\ t \geq 1}$ can be guaranteed i.e., there exists a constant $L_f > 0$ s.t. $\forall\ i \in \small[N \small]$,
    \begin{equation*}\label{eqn:e17}
    \lVert f_{i,t}(p^g_{i,t}) - f_{i,t}(q^g_{i,t})  \rVert \leq L_f \lVert p^g_{i,t} - q^g_{i,t} \rVert,\ \forall\ p^g_{i,t}, q^g_{i,t} \in P_i
    \end{equation*}
    \begin{equation}\label{eqn:e18}
    =>\ \lVert \nabla f_{i,t}(p^g_{i,t})  \rVert \leq L_f,\ \forall\ p^g_{i,t} \in P_i
    \end{equation}
\end{corollary}

\begin{bound} \label{bnd:b1}
Let Assumptions \ref{ass:a1} and \ref{ass:a2} hold. 
Then, iterates $\{\nabla_{p^g_{i,t}}\mathcal{L}_t\}_{i \in \small[N \small],\ t \geq 1}$ are bounded, i.e., there exists a constant s.t. for all $t \geq 1$ and $i \in \small[N \small]$,
\begin{equation}\label{eqn:e23}
    \begin{aligned}
        \lVert \nabla_{p^g_{i,t}}\mathcal{L}_t \rVert &\ = \lVert \nabla f_{i,t}(p^{g}_{i,t}) + \tilde{\lambda}_{i,t} \rVert 
        \\&\
        \leq \lVert \nabla f_{i,t}(p^{g}_{i,t}) \rVert + \lVert \tilde{\lambda}_{i,t} \rVert 
        \\&\ 
        \leq L_f + C_{\lambda}
    \end{aligned}
\end{equation}
\end{bound}

\begin{lemma} \label{lem:l2}
Suppose that for the sequences $\{p_{i,t}^g\}_{i \in \small[N \small],\ t \geq 1}$, there exist a positive scalar sequence $\{\delta_t\}_{t \geq 1}$. Then, for any $T \geq 1$,
\begin{equation}\label{eqn:e24}
    \begin{aligned}
        &\overset{T}{\underset{t=1}{\Sigma}} \overset{N}{\underset{i=1}{\Sigma}} \frac{1}{2 \delta_t} \big[\lVert p_{i,t}^g - p_i^g \rVert^2 - \lVert p_{i,t+1}^g - p_i^g \rVert^2 \big] 
        \\&\ 
        \overset{(a)}{\leq} \frac{1}{2 \delta_1} \overset{N}{\underset{i=1}{\Sigma}} \lVert p_{i,1}^g - p_i^g \rVert^2 + \frac{1}{2} \overset{T}{\underset{t=2}{\Sigma}} \Big(\frac{1}{\delta_t} - \frac{1}{\delta_{t-1}} \Big) \overset{N}{\underset{i=1}{\Sigma}} \lVert p_{i,t}^g - p_i^g \rVert^2 
        \\&\ 
        \overset{(b)}{\leq} \frac{4C_P^2N}{2 \delta_1} + \frac{4C_P^2N}{2} \Big(\frac{1}{\delta_T} - \frac{1}{\delta_1} \Big) = \frac{2NC_P^2}{\delta_T}
    \end{aligned}
\end{equation}
where, $(a)$ implies a negative term $\frac{-1}{2 \delta_T} \overset{N}{\underset{i=1}{\Sigma}} \lVert p_{i,T+1}^g - p_i^g \rVert^2$ being dropped, and $(b)$ is implied from Assumption \ref{ass:a2}. 
\end{lemma}

\begin{corollary} \label{cor:c3}
Following from Lemma \ref{lem:l2}, a stepsize choice of $\delta_t = \frac{1}{\sqrt{t}}$ will yield,
\begin{equation}\label{eqn:e25}
    \overset{T}{\underset{t=1}{\Sigma}} \delta_t = \overset{T}{\underset{t=1}{\Sigma}} \frac{1}{\sqrt{t}} \overset{(a)}{\leq} 1 + \int_1^T \frac{1}{\sqrt{t}} \mathrm{d}x = 2\sqrt{T} - 1
\end{equation}
where (a) comes from \textit{Euler–Maclaurin formula}.
\end{corollary}

\begin{corollary} \label{cor:c4}
Consider sequences $\{\theta_{i,t}\}_{i \in \small[N \small],\ t \geq 1}$ generated by online distributed algorithm step \eqref{eqn:e8}. Then, it can be shown that,
\begin{equation}\label{eqn:e26}
    \begin{aligned}
        \overset{N}{\underset{i=1}{\Sigma}} \overset{N}{\underset{j=1}{\Sigma}} B_{ij} \small[\theta_{j,t} - \theta_j \small] - \overset{N}{\underset{i=1}{\Sigma}} \small[\theta_{i,t} - \theta_i \small] \overset{N}{\underset{j=1}{\Sigma}} B_{ij} = 0
    \end{aligned}
\end{equation}
\end{corollary}

\begin{lemma} \label{lem:l3}
Consider sequences $\{\tilde{\lambda}_{i,t}\}_{i \in \small[N \small],\ t \geq 1}$ and $\{\lambda_{i,t}\}_{i \in \small[N \small],\ t \geq 1}$ generated by online distributed algorithm \eqref{eqn:e9}-\eqref{eqn:e10}. Then, the following chain of relations can be obtained, 
\begin{equation}\label{eqn:e27}
    \begin{aligned}
        \overset{N}{\underset{i=1}{\Sigma}} \lVert \tilde{\lambda}_{i,t} - \bar{\lambda}_t \rVert & \overset{(a)}{\leq} \overset{N}{\underset{i=1}{\Sigma}} \overset{N}{\underset{j=1}{\Sigma}} W_{ij} \lVert \lambda_{j,t} - \bar{\lambda}_t \rVert 
        \\&\
        \overset{(b)}{\leq} \overset{N}{\underset{j=1}{\Sigma}} \lVert \lambda_{j,t} - \bar{\lambda}_t \rVert \overset{N}{\underset{i=1}{\Sigma}} W_{ij} \overset{(c)}{=} \overset{N}{\underset{i=1}{\Sigma}} \lVert \lambda_{i,t} - \bar{\lambda}_t \rVert
    \end{aligned}
\end{equation}
where, (a) follows from the definition of $\tilde{\lambda}_{i,t}$ in \eqref{eqn:e9}, the convexity of the norm and double stochasticity of $W$, while (b) implies reordering of summations and (c) draws upon Assumption \eqref{ass:a4}.
\end{lemma}

\begin{corollary} \label{cor:c5}
Following Lemma \ref{lem:l2}, the aforementioned chain of relations can be established, 
\begin{equation}\label{eqn:e28}
\overset{N}{\underset{i=1}{\Sigma}} \tilde{\lambda}_{i,t} \overset{(a)}{=} \overset{N}{\underset{i=1}{\Sigma}} \overset{N}{\underset{j=1}{\Sigma}} W_{ij} \lambda_{j,t} \overset{(b)}{=} \overset{N}{\underset{j=1}{\Sigma}} \lambda_{j,t} \overset{N}{\underset{i=1}{\Sigma}} W_{ij} \overset{(c)}{=} \overset{N}{\underset{i=1}{\Sigma}} \lambda_{i,t}
\end{equation}
where (a) comes from \eqref{eqn:e9}, (b) is due to the reordering of summation terms and (c) implies Assumption \ref{ass:a4}.
\end{corollary}

\begin{bound} \label{bnd:b2}
Consider sequence $\{\tilde{\lambda}_{i,t}\}_{i \in \small[N \small],\ t \geq 1}$ generated by online distributed algorithm \eqref{eqn:e9}. Then, for any $t \geq 1$ and $i \in \small[N \small]$, 
\begin{equation}\label{eqn:e29}
    \lVert \tilde{\lambda}_{i,t} \rVert = \lVert \underset{j \in \mathcal{N}_i}{\Sigma} W_{ij} \lambda_{j,t} \rVert \overset{(a)}{\leq} C_{\lambda} \overset{N}{\underset{j=1}{\Sigma}} W_{ij} = C_{\lambda}
\end{equation}
where (a) follows from Assumption \ref{ass:a3}-\ref{ass:a4}.
\end{bound}

\begin{bound} \label{bnd:b3}
Let Assumption \ref{ass:a3} hold. Then, iterates $\{\nabla_{\theta_{i,t}}\mathcal{L}_t\}_{i \in \small[N \small],\ t \geq 1}$ are bounded, i.e., there exists a constant s.t. for all $t \geq 1$ and $i \in \small[N \small]$, 
\begin{equation}\label{eqn:e30}
    \begin{aligned}
        \lVert \nabla_{\theta_{i,t}}\mathcal{L}_t \rVert = \lVert -\tilde{\lambda}_{i,t}\underset{j \in \mathcal{N}_i}{\Sigma}B_{ij} + \underset{j \in \mathcal{N}_i}{\Sigma}B_{ij}\tilde{\lambda}_{j,t} \rVert \leq 2C_{\lambda} \lVert B \rVert_{\infty}
    \end{aligned}
\end{equation}
\end{bound}

\begin{bound} \label{bnd:b4}
Let Assumption \ref{ass:a2} hold. Then, iterates $\{\nabla_{\lambda_{i,t}}\mathcal{L}_t\}_{i \in \small[N \small],\ t \geq 1}$ are bounded, i.e., there exists a constant s.t. for all $t \geq 1$ and $i \in \small[N \small]$,
\begin{equation}\label{eqn:e31}
    \begin{aligned}
        \lVert \nabla_{\lambda_{i,t}}\mathcal{L}_t \rVert &\ = \lVert p^g_{i,t} - p^l_{i,t} - \theta_{i,t}\underset{j \in \mathcal{N}_i}{\Sigma}B_{ij} + \underset{j \in \mathcal{N}_i}{\Sigma}B_{ij} \theta_{j,t} \rVert 
        \\&\
        \leq C_P + 2C_{\theta} \lVert B \rVert_{\infty}
    \end{aligned}
\end{equation}
\end{bound}

\begin{lemma} \label{lem:l4}
Suppose that for the sequences $\{\lambda_{i,t}\}_{i \in \small[N \small],\ t \geq 1}$, there exist a positive scalar sequence $\{\alpha_t\}_{t \geq 1}$. Then, for any $T \geq 1$,
\begin{equation}\label{eqn:e32}
    \begin{aligned}
        &\ \overset{T}{\underset{t=1}{\Sigma}} \overset{N}{\underset{i=1}{\Sigma}} \frac{1}{2 \alpha_t} \big[\lVert \lambda_{i,t} - \lambda \rVert^2 - \lVert \lambda_{i,t+1} - \lambda \rVert^2 \big]
        \\&\ 
        \overset{(a)}{\leq} \frac{1}{2 \alpha_1} \overset{N}{\underset{i=1}{\Sigma}} \lVert \lambda_{i,1} - \lambda \rVert^2 + \frac{1}{2} \overset{T}{\underset{t=2}{\Sigma}} \Big(\frac{1}{\alpha_t} - \frac{1}{\alpha_{t-1}} \Big) \overset{N}{\underset{i=1}{\Sigma}} \lVert \lambda_{i,t} - \lambda \rVert^2 
        \\&\ 
        \overset{(b)}{\leq} \frac{4C_{\lambda}^2N}{2 \alpha_1} + \frac{4C_{\lambda}^2N}{2} \Big(\frac{1}{\alpha_T} - \frac{1}{\alpha_1} \Big) = \frac{2NC_{\lambda}^2}{\alpha_T}
    \end{aligned}
\end{equation}
where, $(a)$ implies a negative term $\frac{-1}{2 \alpha_T} \overset{N}{\underset{i=1}{\Sigma}} \lVert \lambda_{i,T+1} - \lambda \rVert^2$ being dropped, and $(b)$ is implied from Assumption \ref{ass:a3}.
\end{lemma}

\begin{corollary} \label{cor:c6}
Following from \eqref{eqn:e10}, a stepsize choice of $\beta_t = \frac{1}{t}$ will yield,
\begin{equation}\label{eqn:e33}
    \overset{T}{\underset{t=1}{\Sigma}} \beta_t = \overset{T}{\underset{t=1}{\Sigma}} \frac{1}{t} \overset{(a)}{\leq} 1 + \int_1^T \frac{1}{t} \mathrm{d}x = 1+ ln(T)
\end{equation}
where (a) comes from \textit{Euler's Constant approximation}.
\end{corollary}

\begin{corollary} \label{cor:c7}
Following from \eqref{eqn:e42}, in order to achieve a sublinear regret, a stepsize choice of $\alpha_t = \frac{1}{\sqrt{t}}$ and $\beta_t = \frac{1}{t}$, would yield the following relations :
\begin{equation}\label{eqn:e33b}
    \overset{T}{\underset{t=1}{\Sigma}} \frac{\beta_t}{\alpha_t} = \overset{T}{\underset{t=1}{\Sigma}} \frac{1}{\sqrt{t}} \overset{(a)}{\leq} 1 + \int_1^T \frac{1}{\sqrt{t}} \mathrm{d}x = 2\sqrt{T} - 1
\end{equation}
\begin{equation}\label{eqn:e33c}
    \overset{T}{\underset{t=1}{\Sigma}} \frac{\beta_t^2}{\alpha_t} = \overset{T}{\underset{t=1}{\Sigma}} \frac{1}{t \sqrt{t}} \overset{(a)}{\leq} 1 + \int_1^T \frac{1}{t \sqrt{t}} \mathrm{d}x = 3 - \frac{2}{\sqrt{T}}
\end{equation}
where (a) comes from \textit{Euler–Maclaurin formula}.
\end{corollary}

\begin{proposition} \label{prop:p1}
\cite[Lemma 5.2]{sundhar2012new} Consider a graph $\mathcal{G}(\mathcal{V},\mathcal{E},W)$ satisfying Assumption \ref{ass:a4} and a set of sequences $x_{i,t}\ \forall\ i \in \small[N \small]$ defined by - 
\begin{equation}\label{eqn:e34}
    x_{i,t+1} = \overset{N}{\underset{i=1}{\Sigma}} W_{ij} x_{j,t} + \epsilon_{i,t+1} \ \forall\ t \geq 1
\end{equation}
Let $\bar{x}_t$ denote the average of $x_{i,t}$ for $i \in \small[N \small]$, i.e., $\bar{x}_t$ = $\frac{1}{N}\overset{N}{\underset{i=1}{\Sigma}}x_{i,t}$. Then,
\begin{equation}\label{eqn:e35}
    \begin{aligned}
        \lVert x_{i,t+1} - \bar{x}_{t+1} \rVert\ \leq\ & N\kappa \omega^t \underset{j}{max} \lVert x_{j,1} \rVert + \kappa \overset{t-1}{\underset{l=1}{\Sigma}} \omega^{t-l} \overset{N}{\underset{j=1}{\Sigma}} \lVert \epsilon_{j,l+1} \rVert
        \\&\ 
        + \frac{1}{N} \overset{N}{\underset{j=1}{\Sigma}} \lVert \epsilon_{j,t+1} \rVert + \lVert \epsilon_{i,t+1} \rVert
    \end{aligned}
\end{equation}
where, $\kappa = (1-\frac{\eta}{2N^2})^{-2}$ and, $\omega = (1-\frac{\eta}{2N^2})$.\\
The above lemma is modified for the case of  time-invariant undirected graph for our analysis.
\end{proposition}  

\begin{corollary} \label{cor:c8}
Following Proposition \ref{prop:p1}, the iterate $\{\lambda_{i,t}\}_{i \in \small[N \small],\ t \geq 1}$ in algorithm \eqref{eqn:e9}-\eqref{eqn:e10}, can be represented by the relation \eqref{eqn:e34} and can be rewritten as, 
\begin{equation}\label{eqn:e36}
    \lambda_{i,t+1} = \overset{N}{\underset{i=1}{\Sigma}} W_{ij} \lambda_{j,t} + \epsilon_{i,t+1}
\end{equation}
with, $\epsilon_{i,t+1} = \alpha_t \nabla_{\lambda_{i,t}}\mathcal{L}_t - \beta_t \nabla_{\theta_{i,t}}\mathcal{L}_t$.\\
Suppose, that for the sequences $\{\lambda_{i,t}\}$ defined above, there exist positive scalar sequences $\{\alpha_t\}_{t \geq 1}$ and $\{\beta_t\}_{t \geq 1}$ s.t. 
\begin{equation}\label{eqn:e37}
    \lVert \epsilon_{i,t} \rVert \leq (C_P + 2C_{\theta} \lVert B \rVert_{\infty}) \alpha_t + (2 C_{\lambda} \lVert B \rVert_{\infty}) \beta_t,\ \forall\ i \in \small[N \small],\ t \geq 1 
\end{equation}
Then, for any $T \geq 2$,
\begin{equation}\label{eqn:e38}
    \begin{aligned}
        \overset{T-1}{\underset{t=1}{\Sigma}}& \lVert \lambda_{i,t+1} - \bar{\lambda}_{t+1} \rVert
        \\&\ 
        \leq \kappa N \underset{j}{max} \lVert \lambda_{j,1} \rVert \overset{T-1}{\underset{t=1}{\Sigma}} \omega^t + \kappa \overset{T-1}{\underset{t=1}{\Sigma}} \overset{t-1}{\underset{l=1}{\Sigma}} \omega^{t-l} \overset{N}{\underset{j=1}{\Sigma}} \lVert \epsilon_{j,l+1} \rVert \\&\ + \frac{1}{N} \overset{T-1}{\underset{t=1}{\Sigma}} \overset{N}{\underset{j=1}{\Sigma}} \lVert \epsilon_{j,t+1} \rVert + \overset{T-1}{\underset{t=1}{\Sigma}} \lVert \epsilon_{i,t+1} \rVert 
        \\&\
        \overset{(a)}{\leq} \frac{\kappa N \omega C_{\lambda}}{1 - \omega} + 2(C_P + 2C_{\theta} \lVert B \rVert_{\infty}) \overset{T-1}{\underset{t=1}{\Sigma}} \alpha_t\ + 
        \\&\
        4C_{\lambda} \lVert B \rVert_{\infty} \overset{T-1}{\underset{t=1}{\Sigma}} \beta_t + \kappa N \overset{T-1}{\underset{t=1}{\Sigma}} \overset{t-1}{\underset{l=1}{\Sigma}} \omega^{t-l} \big[(C_P + 2C_{\theta} \lVert B \rVert_{\infty})\alpha_t
        \\&\
        + (2C_{\lambda}\lVert B \rVert_{\infty})\beta_t\big] 
        \\&\
        \overset{(b)}{=} \frac{\kappa N \omega C_{\lambda}}{1 - \omega} + 2(C_P + 2C_{\theta} \lVert B \rVert_{\infty}) \overset{T-1}{\underset{t=1}{\Sigma}} \alpha_t\ + 
        \\&\
        4C_{\lambda} \lVert B \rVert_{\infty} \overset{T-1}{\underset{t=1}{\Sigma}} \beta_t + \kappa N (C_P + 2C_{\theta} \lVert B \rVert_{\infty}) \overset{T-2}{\underset{l=1}{\Sigma}} \omega^l \overset{T-l-1}{\underset{t=1}{\Sigma}} \alpha_t 
        \\&\
        + (2\kappa N C_{\lambda}\lVert B \rVert_{\infty})\overset{T-2}{\underset{l=1}{\Sigma}} \omega^l \overset{T-l-1}{\underset{t=1}{\Sigma}} \beta_t 
        \\&\
        \overset{(c)}{\leq} \frac{\kappa N \omega C_{\lambda}}{1 - \omega} + 2(C_P + 2C_{\theta} \lVert B \rVert_{\infty}) \overset{T-1}{\underset{t=1}{\Sigma}} \alpha_t\ + 
        \\&\
        4C_{\lambda} \lVert B \rVert_{\infty} \overset{T-1}{\underset{t=1}{\Sigma}} \beta_t + \kappa N (C_P + 2C_{\theta} \lVert B \rVert_{\infty}) \overset{T-2}{\underset{l=1}{\Sigma}} \omega^l \overset{T-1}{\underset{t=1}{\Sigma}} \alpha_t
        \\&\
        + (2\kappa N C_{\lambda}\lVert B \rVert_{\infty})\overset{T-2}{\underset{l=1}{\Sigma}} \omega^l \overset{T-1}{\underset{t=1}{\Sigma}} \beta_t
    \end{aligned}
\end{equation}
where, (a) follows from infinite geometric sum bound and \eqref{eqn:e36}, (b) involves the reordering of summation terms and (c) implies positive terms i.e. $\kappa N (C_P + 2C_{\theta} \lVert B \rVert_{\infty}) \overset{T-2}{\underset{l=1}{\Sigma}} \omega^l \overset{T-1}{\underset{t=T-l}{\Sigma}} \alpha_t + (2\kappa N C_{\lambda}\lVert B \rVert_{\infty})\overset{T-2}{\underset{l=1}{\Sigma}} \omega^l \overset{T-1}{\underset{t=T-l}{\Sigma}} \beta_t$, being added. 
\end{corollary}

\begin{corollary} \label{cor:c9}
Following Corollary \ref{cor:c6}, we can establish the following relation,
\begin{equation}\label{eqn:e39}
    \begin{aligned}
        \overset{T}{\underset{t=1}{\Sigma}} \overset{N}{\underset{i=1}{\Sigma}} \lVert \tilde{\lambda}_{i,t} - \bar{\lambda}_{t} \rVert &\ \overset{(a)}{\leq} \overset{T}{\underset{t=1}{\Sigma}} \overset{N}{\underset{i=1}{\Sigma}} \lVert \lambda_{i,t} - \bar{\lambda}_{t} \rVert 
        \\&\
        = \overset{N}{\underset{i=1}{\Sigma}} \lVert \lambda_{i,1} - \bar{\lambda}_{1} \rVert + \overset{T-1}{\underset{t=1}{\Sigma}} \overset{N}{\underset{i=1}{\Sigma}} \lVert \lambda_{i,t+1} - \bar{\lambda}_{t+1} \rVert 
        \\&\ 
        \overset{(b)}{\leq} 2NC_{\lambda} + \frac{\kappa N^2 \omega C_{\lambda}}{1-\omega} + 4NC_{\lambda} \lVert B \rVert_{\infty} \overset{T}{\underset{t=1}{\Sigma}} \beta_t
        \\&\
        + 2N(C_P + 2C_{\theta} \lVert B \rVert_{\infty}) \overset{T}{\underset{t=1}{\Sigma}} \alpha_t\ +  
        \\&\
        \kappa N^2 (C_P + 2C_{\theta} \lVert B \rVert_{\infty}) \overset{T-2}{\underset{l=1}{\Sigma}} \omega^l \overset{T}{\underset{t=1}{\Sigma}} \alpha_t\ + 
        \\&\
        (2\kappa N^2 C_{\lambda}\lVert B \rVert_{\infty})\overset{T-2}{\underset{l=1}{\Sigma}} \omega^l \overset{T}{\underset{t=1}{\Sigma}} \beta_t
    \end{aligned}
\end{equation}
where, (a) is implied from \eqref{eqn:e27} and (b) follows from \eqref{eqn:e38}.
\end{corollary}

\bibliographystyle{ieeetr}
\bibliography{Bibliography.bib}
\balance
\end{document}